\documentclass[reqno,10pt]{amsart}
\usepackage{amsmath,amsxtra,amssymb,latexsym, amscd,amsthm}
\usepackage[unicode]{hyperref}
\usepackage{array,tabularx,longtable,multicol,indentfirst,fancybox,color}
\usepackage{graphicx}
\usepackage{multicol}
\usepackage{mathrsfs}
\usepackage{enumerate}
\advance\hoffset-1truecm\relax

\newtheorem{proposition}{Proposition}[section]
\newtheorem{theorem}{Theorem}[section]
\newtheorem{definition}{Definition}
\newtheorem{lemma}{Lemma}[section]

\newtheorem{remark}{Remark}
\newtheorem{counterexample}{Counterexample}
\newtheorem{example}{Example}[section]
\numberwithin{equation}{section}
\newtheoremstyle{named}{}{}{\itshape}{}{\bfseries}{.}{.5em}{\thmnote{#3's }#1}
\theoremstyle{named}
\newtheorem*{theoremA}{Theorem A}
\newtheorem*{theoremB}{Theorem B}

\def\O{\Omega}
\def\wed{\wedge}
\def\n{\noindent}
\def\C^n{\mathbb C^n}
\def\var{\varepsilon}

\def\wide{\widetilde}
\def\varp{\varphi}

\numberwithin{equation}{section}
\begin{document}
	\title[On Lelong numbers of plurisubharmonic functions on complex spaces]{On Lelong numbers of plurisubharmonic functions on\\ complex spaces}
	\date{\today}	
	\author[Le Mau Hai, Pham Hoang Hiep and Trinh Tung]{Le Mau Hai$^{1*}$, Pham Hoang Hiep$^{2}$, Trinh Tung$^{3}$}
%	\thanks{*Corresponding author: \textsc{Le Mau Hai}}

\thanks{\textit{\textbf{2020 Mathematics Subject Classification.}} 32C15, 32C25, 32C30, 32U25, 32U40, 32W20.} 
\thanks{\textit{\textbf{Keywords.}} 
	 Strong locally  irreducible complex spaces, Regular points, Singular points, Irreducible complex spaces, Locally irreducible complex spaces,  Plurisubharmonic functions on complex spaces, Weakly plurisubharmonic functions on complex spaces, Lelong number of a plurisubharmonic function.}
	\maketitle	
	
\begin{center}
\scriptsize	
	$^1$Department of Mathematics, Hanoi National University of Education,\\ 136-Xuan Thuy Road, Hanoi, Vietnam.\\
	\textit{E-mail:} \texttt{mauhai@hnue.edu.vn}\\
	$^2$ICRTM, Institute of Mathematics, Vietnam Academy of Science and Technology,\\ 18-Hoang Quoc Viet Road, Hanoi, Vietnam.\\ \textit{E-mail:} \texttt{phhiep@math.ac.vn}\\
	$^3$Department of Mathematics, Faculty of Sciences,  Electric Power University,\\ 235-Hoang Quoc Viet Road, Hanoi, Vietnam.\\ \textit{E-mail:} \texttt{tungtrinhvn@gmail.com}
\end{center}

\begin{abstract}
	In this paper, we introduce the notion of strong locally irreducible complex spaces $\widetilde{X}$. Based on this notion we prove the equality $\bar{\nu}_{\varphi}(x)=$ mult$(\widetilde{X},x). \nu_{\varphi}(x)$ for all $x\in \widetilde{X}$, where $\bar{\nu}_{\varphi}(x)$ is the projective mass of a plurisubharmonic function $\varphi$ at $x$ and mult$(\widetilde{X},x)$ is the multiplicity of $\widetilde{X}$ at $x$ and $\nu_{\varphi}(x)$ is Lelong number of $\varphi$ at $x$. Moreover, we show that the closure of the upper-level sets $\{z\in \widetilde{X}:\nu_{\varphi}(z)\geq c\}$ of a plurisubharmonic function $\varphi$ on a strong locally irreducible complex space $\widetilde{X}$ is a subvariety of $\widetilde{X}$ for all $c\geq 0$.

	%\noindent \textsc{Résumé.} Dans cet article, nous introduisons la notion d'espaces complexes forts localement irréductibles $\widetilde{X}$. En nous appuyant sur cette notion, nous prouvons l'égalité $\bar{\nu}_{\varphi}(x)=\textup{mult}(\wide{X},x).\nu_{\varphi}(x)$ pour all $x\in \wide{X}$, où $\bar{\nu}_{\varphi}(x)$ est la masse projective de la fonction plurisubharmonique $\varphi$ à $x$ et $\textup{mult}(\wide{X},x)$ est la multiplicité de $\wide{X}$ à $x$ et $\nu_{\varphi}(x)$ est le nombre de Lelong de $\varphi$ à $x $. De plus, nous montrons que la fermeture des ensembles de niveau supérieur $\{z\in \wide{X}:\nu_{\varp}(z)\geq c\}$ d'une fonction plurisubharmonique $\varp$ sur un l'espace complexe fort localement irréductible $\wide{X}$ est une sous-variété de $\wide{X}$ pour tout $c\geq 0$.
\end{abstract}
%\tableofcontents

\section{Introduction}
A complex space $X$ is a topological space with an atlas of charts which are isomorphic to  analytic sets in complex Euclidean spaces such that the transition maps are holomorphic in the sense of holomorphic maps between analytic sets. A point $z$ of a complex space $X$ is said to be regular if there is a neighborhood $U\ni z$ in $X$ such that $U$ is isomorphic to an open subset of a complex Euclidean space. The set of regular points of $X$ is denoted by $X_{reg}$. We define the set of singularity points of $X$ by $X_{sing}:=X\backslash X_{reg}$. A complex space $X$ is called to be reducible if $X=X_1\cup X_2$, where $X_1, X_2$ are proper subvarieties of $X$. If $X$ can not be represented in this form, $X$ is said to be irreducible. Of course, $X$ is irreducible if and only if $X_{reg}$ is connected. A complex space $X$ has a decomposition $X=\bigcup\limits_{i\in I} X_i$, where $X_i$ are irreducible complex spaces, $X_i\not\subset X_j$ for all $i,j\in I$, $i\not = j$ and the family $\{ X_i \}_{i\in I}$ is locally finite. Let $X,Y$ be complex spaces. A continuous map $f:X\to Y$ is said to be a holomorphic map if for each point $a\in X$ there exists a neighborhood $U$ of $a$ and a neighborhood $V$ of $f(a)$ such that $f:U\to V$ is a holomorphic map in the sense of a holomorphic map between analytic sets. A function $\varphi : X\to [-\infty, +\infty )$ is said to be weakly plurisubharmonic on $X$ if $\varphi$ plurisubharmonic on $X_{reg}$ and for all $x_0\in X$
$$\varlimsup_{ x\in X_{reg}, x\to x_0}\varphi(x) = \varphi(x_0),$$  
where $\varphi$ is said to be plurisubharmonic on $X_{reg}$ if for all $a\in X_{reg}$ there exists a neighbourhood $U$ of $a$ in $X$ such that $\psi:U\cong V$ is an analytically isomorphic to an open subset $V\subset \C^n$ and there exists a plurisubharmonic function $g:V\rightarrow [-\infty,+\infty)$ with $\varphi= g\circ\psi$.
The set of weakly plurisubharmonic functions (resp. negative weakly plurisubharmonic functions) on $X$ is denoted by $PSH_{w} (X)$ (resp. $PSH_{w}^- (X)$). 

A function $\varphi : X\to [-\infty, +\infty )$ is said to be plurisubharmonic on $X$ if $\varphi$ is upper semi-continuous on $X$ and $\varphi\circ f$ is subharmonic on the unit disk $\triangle=\{z\in\mathbb{C}: |z|<1\}$, for all holomorphic mapping $f:\triangle\rightarrow X$. The set of plurisubharmonic functions (resp. negative plurisubharmonic functions) on $X$ is denoted by $PSH (X)$ (resp. $PSH^- (X)$). Note that the above definition coincides with the notion of weakly plurisubharmonic functions on a complex space $X$ introduced in \cite{FN80}. The purpose of this paper is to study some properties of Lelong number (slope) of a plurisubharmonic function $\varphi$ at a point $x\in X$, denoted by $\nu_{\varphi}(x)$,  on a complex space $X$ introduced and investigated by Demailly in \cite{De85}. To get the main results in this paper, we need to add new notions and new techniques. Namely, in Section \ref{section3}, for a plurisubharmonic function $\varphi$ on a complex space $X$, we associate it with the two new plurisubharmonic functions $\varphi_{aver}$ and $\varphi_{max}$ by using ramified coverings described in \cite{De12} (see the local parametrization theorem (II.4.9) in \cite{De12}). Using these auxiliary plurisubharmonic functions, we can apply well-established techniques of plurisubharmonic functions on open subsets of $\C^n$ and thanks to that, it helps with studying results related to plurisubharmonic functions on complex spaces with singularities. Next, in Section \ref{section2}, we introduce the notion about strong locally  irreducible complex spaces $\wide{X}$ (detail, see Definition \ref{def19} in Section \ref{section2}). With new plurisubharmonic functions $\varphi_{aver}, \varphi_{max}$ and the strong locally irreducible complex space $\wide{X}$, we may state main results of this paper as follows. The first main result (Theorem \ref{theo38}) aims to describe the relation between the Lelong numbers of three plurisubharmonic functions $\varphi, \varphi_{aver}, \varphi_{max}$ on a strong locally  irreducible complex space $\wide{X}$.
\begin{theoremA}\label{theoA}
	Let $X$ be a complex space and $\varphi: X\rightarrow [-\infty,+\infty)$ be a  plurisubharmonic function on $X$. Then for every $x\in X$ the following conclusions holds
	\vskip0.1cm
	{\bf i)} $\nu_{\varphi_{aver}}(x)\geq \nu_{\varphi}(x) = \nu_{  \varphi_{max}}(x).$
	\vskip0.1cm
	{\bf ii)} If $\wide{X}$ is a strong locally irreducible complex space then for all $x\in\wide{X}$ the following equalities holds
	\begin{equation*}\label{eq1.2}
	\nu_{\varphi_{aver}}(x) = \nu_{\varphi}(x) = \nu_{  \varphi_{max}}(x).
	\end{equation*}
\end{theoremA}
\n It is worth noting that without the hypothesis about the strong local irreducibility of $X$, then the first equality in {\bf i)} is not true. We have given a counter-example to clarify this claim (Counterexample \ref{phvd1}).

To state the second main result, we recall the notion about Lelong number of a plurisubharmonic function $\varphi$ on a complex space $X$ and its projective mass. Let $X$ be a complex space and $\varphi$ be a plurisubharmonic function (abbreviated as \textit{psh function}) defined on $X$. The Lelong number (slope) of $\varphi$ at $z\in X$ is defined as the slope with respect to log at $-\infty$ and given by
$\nu_{\varphi}(z):=\lim\inf\limits_{x\to z}\frac{\varphi(x)}{\log \|f(x)-f(z)\|}$, where $f: U\to A$ is a chart from a neighborhood $U$ of $z$ to an analytic set $A$ in some complex Euclidean space $\mathbb{C}^m$ 
(\cite{De85}, p. 45). Another definition is provided by the projective mass of a plurisubharmonic function $\varphi$ defined on a complex space $X$ at a point $x\in X$ introduced and investigated by Demailly in (\cite{De82}, Déf. 3) and given by the following formula 
$$\bar{\nu}_{\varphi }(x)=\lim\limits_{r\to 0}\int\limits_{\|z-x\|<r} dd^c\varphi\wedge(dd^c\log\|z-x\|)^{k-1}\wedge[X],$$

\n where $k=\dim_{\mathbb{C}}X$ is the dimension of complex space $X$. These two definitions are equivalent for psh functions on a smooth domain, but in the singular setting, the relation between them is less clear. In the case $X$ is a complex space, by Demailly's comparison theorem (see \cite{De82}, Th\'eor\`eme 4), we note that the following inequality holds
\begin{equation}\label{eq1.4}
\bar{\nu}_{\varphi}(x)\geq \textup{mult}(X,x).\nu_{\varphi}(x).
\end{equation}

\n for all $x\in X$, where $\textup{mult}(X,x)$ is the multiplicity of $X$ at $x$ (Also see Remark A.5 in \cite{BBEGZ19}).

The second main result in this paper (Theorem \ref{theo311}) is if $X=A$ is a $k$-dimensional analytic subset in $\C^n$  and $\varphi\in PSH(A)$ is a psh function on $A$. Consider the function $\varphi_{aver}$ associated with $\varphi$. Then we obtain the equality in \eqref{eq1.4}. Namely, we have the following result
\begin{theoremB}\label{theoB}
	Let $A$ be a $k$-dimensional analytic subset in $\C^n, n\geq k$ and $\varphi\in PSH(A)$ be a psh function on $A$. Then the following equality holds for all $a\in A$
	$$\bar{\nu}_{\varphi}(a)=\textup{mult}(A,a).\nu_{\varphi_{aver}}(a).$$  
	\n Furthermore, if $\wide{A}$ is a $k$-dimensional strong  locally irreducible  analytic subset of $\C^n$, then for all $a\in \wide{A}$ 
	$$\bar{\nu}_{\varphi}(a)=\textup{mult}(\wide{A},a).\nu_{\varphi}(a).$$
\end{theoremB}
\noindent Notice that, Pan's recent result (see Main Theorem in \cite{Pan2025}) showed that if X is a locally irreducible complex space, then $\bar{\nu}_{\varphi}(x)\leq C_{x}. \textup{mult}(X,x).\nu_{\varphi}(x)$, for every $x \in X$. By the above \textbf{Theorem B}, under the additional assumption that X is a strong locally irreducible complex space we achieved the inequality sign turns into an equality.

Finally, we prove a result of a theorem of Siu's type which holds on strong locally irreducible complex spaces. This is Proposition \ref{pro312} asserting that the closure of the upper-level sets $\{z\in X: \nu_{\varphi}(z)\geq c\}$ of a psh function $\varphi$ defined on a  strong locally irreducible complex space $X$ is a subvariety of $X$ for all $c\geq 0$.

This article is divided into four sections and organized as follows. In Section \ref{section2} we recall the definition of complex spaces and analytic subsets in complex spaces and some its basic results. Among notions concerning to analytic subsets in complex spaces we deal with regular, singular points on complex spaces. At the same time, we recall the set of irreducible points and reducible points and holomorphic mappings between complex spaces. The important definition of this section is Definition \ref{def19} about strong locally irreducible complex spaces $\wide{X}$. Section \ref{section3} is devoted to the presentation of plurisubharmonic functions on complex spaces. Here we are interested in  the two definitions of plurisubharmonic functions. These are weakly plurisubharmonic functions and plurisubharmonic functions. Besides, we introduce the two plurisubharmonic functions $\varphi_{aver}$ and $\varphi_{max}$ associated with a weakly plurisubharmonic function $\varphi$ on a complex space $X$. In Section \ref{section4} we recall the notion of the  Lelong number of a plurisubharmonic function on a complex space $X$ introduced and  investigated by Demailly in \cite{De85}. We establish some basic properties of this quantity. In this section we prove the two essential results mentioned as above. Finally, at Proposition \ref{pro312},  we show that the closure of the upper-level set $\{z \in X:\varphi(z)\geq c\}$ of a plurisubharmonic function $\varphi$ on a strong locally irreducible complex space $X$ is a subvariety of $X$ for all $c>0$.   

\section{Some notations and basic results of analytic subsets on complex spaces}\label{section2}

\n Some elements of the theory of complex manifolds, complex spaces, analytic subsets on manifolds and on complex spaces, holomorphic mappings between complex spaces we refer readers to \cite{Chir89,De12,Ge76,GrRe84,GuRo}. Elements of theory of complex analysis of several variables, $\overline{\partial}$-equation, theory of Stein manifolds and its applications we refer readers to the monograph of L. H\"ormander in \cite{Ho90}. Elements of pluripotential theory on open subsets of $\C^n$, complex Monge-Amp\`ere operator, complex Monge-Amp\`ere equations and related problems we refer readers to \cite{ACKPZ09,BT76,BT82,CKNS85,Ce98,Ce04,De92,De93,Ho94,Kl91,Ko98,Ko05}. For the notion of Lelong number and related problems readers can find in \cite{De12,Lelong57,Kis94,Sk72,Siu74}. For the notions of currents in complex analysis, positive currents and extended problems concerning with currents we refer readers to paper: \cite{DiSi06,DiSi07,DiSi09,DiSi10}. Besides, readers also can to reach some notions about plurisubharmonic functions on complex spaces in \cite{FN80}.

The following definition of  complex manifolds is inspired from  \cite{GriHar78}.

\begin{definition}\label{de11} Let $X$ be a topological manifold. $X$ is said to be a complex manifold if there exists an atlas of charts $(\varphi_{\alpha}, U_{\alpha})_{\alpha\in I}$ such that $\{U_{\alpha}\}_{\alpha\in I}$ are an open covering of $X$ and  $\varphi_{\alpha}: U_{\alpha}\cong \varphi_{\alpha}(U_{\alpha})\subset\C^n$ are isomorphic to open sets $\varphi_{\alpha}(U_{\alpha})$ in $\mathbb{C}^n$ such that the transition maps $\varphi_{\alpha}\circ \varphi^{-1}_{\beta}: \varphi_{\beta}(U_{\alpha}\cap U_{\beta})\rightarrow \varphi_{\alpha}(U_{\alpha}\cap U_{\beta})$ are holomorphic  for all $\alpha, \beta\in I$.
\end{definition}

\n Now we give the notion of  an analytic set in a complex manifold as in \cite{Chir89}.

\begin{definition}\label{de12} Let $X$ be a complex maniflold and $A\subset X$. $A$ is said to be an analytic subset of $X$ if for each point $a\in X$ there exists a neighborhood $U$ of $a$ , and  holomorphic functions $f_1,\ldots,f_m$ on $U$ such that
	$$A\cap U=\{z\in U: f_1(z)=\cdots=f_m(z)=0\}.$$ 
\end{definition}

\n From now until to the end of this paper when we deal with  holomorphic functions on a complex manifolds (resp. on complex spaces), that means that these are $\mathbb{C}$-valued functions which are holomorphic in the sense in \cite{GuRo} or \cite{Ge76,GrRe84}. 

The following results is a consequence of $\bar\partial$-equation with $L^2$ estimates (see \cite{Ho90})

\begin{proposition}\label{pro11}
	Let $\Omega$ be a psedoconvex domain in $\mathbb{C}^n$ and A be an analytic subset of $\O$. Then for every open subset $D\Subset\Omega$ there exist holomorphic functions $f_1,\ldots,f_m$ on $D$ such that
	$$A\cap D=\{z\in D: f_1(z)=\cdots=f_m(z)=0\}.$$ 
\end{proposition}

\n The above result can be extended to the case $X$ is a Stein manifold (see Theorem B Cartan in \cite{GuRo} and \cite{Ho90}).

\begin{proposition}\label{pro12}
	Let $X$ be a Stein manifold and $A$ be an analytic subset of $X$. Then for every open subset $D\Subset X$ there exist holomorphic functions  $f_1,\ldots,f_m$ on $D$ such that
	$$A\cap D=\{z\in D: f_1(z)=\cdots=f_m(z)=0\}.$$ 
\end{proposition}

\n The next definition is important in the sequel.

\begin{definition}\label{de13} A point $a$ of an analytic subset $A$ of a complex manifold $X$ is said to be regular if there is a neighborhood $U\ni a$ in $X$ such that $A\cap U$ is a complex submanifold of $X$. The complex dimension of this submanifold is called to be the dimension of the set $A$ at its regular $a$, and is denoted by $dim_a A$. The set of regular points of an analytic subset $A$ is denoted by $A_{reg}$. Every point in the complement $A\setminus A_{reg}=:A_{sing}$ is called to be a singular point of $A$. 
\end{definition}
\begin{example}\label{vd1}{\rm
		\n In $\mathbb{C}^2$ consider the analytic subset $A=\{(z_1,z_2): z_1 z_2=0\}$. Then every point of the form $a= (z_1,0)$ with $z_1\ne 0$ (resp. $a=(0,z_2)$, $ z_2\ne 0$) are regular points of $A$. The point $(0,0)$ is the singular point of $A$. It is easy to see that if $a\in A$, $a=(z_1,0), z_1\ne 0$ (resp. $a=(0,z_2), z_2\ne 0$) then $dim_{a}A=1$.}
\end{example}

\n We recall the following result which has been proved in \cite{GrRe84}, p. 117.	

\begin{proposition}\label{pro13} 
	Let $X$ be a complex manifold and $A$ be an analytic subset of $X$. Then $A_{sing}$ is an analytic subset of $X$.
\end{proposition}

\n We refer readers to the following definition in \cite{Chir89}, p. 54-55.
\begin{definition}\label{de14}
	Let $A$ be an analytic subset of a complex manifold $X$. $A$ is said to be reducible if $A=A_1\cup A_2$, where $A_1,A_2$ are also analytic subsets of $X$, and $A_1\not = A$, $A_2\not = A$. If $A$ cannot be represented in this form, $A$ is called irreducible in $X$.	
\end{definition}

If $A$ is an analytic subset of $X$ then it has a decomposition $A=\bigcup\limits_{i=1}^k A_i$, where $A_i$ are irreducible analytic subsets, $A_i\nsubseteq A_j$ when $i\ne j$ (they are called to be irreducible branches of $A$). The above decomposition is unique.

\n Moreover, by a result in \cite{Chir89}(see Section 5.3, p. 55) an analytic set $A$ is irreducible if and only if the set $A_{reg}$ is connected.  
\begin{definition}\label{de15}
	An analytic set $A$ is irreducible at $a\in A$ if there exists a basis of neighbourhoods $U_j$ of $a$ such that all analytic sets $A\cap U_j$ are irreducible in $U_j$.	
\end{definition}

\n The following definition can be found in \cite{Chir89} (p. 19).

\begin{definition}\label{de16}
	Let $A$ be an analytic subset of a complex manifold $X$. The dimension of $A$ at an arbitrary point $a\in A$ is defined 
	$$dim_a A=\varlimsup\limits_{z\in A_{reg}, z\to a} dim_z A.$$
	
	\n The dimension of $A$ is, by definition, the maximum of its dimension at points
	$$dim A=\max\limits_{z\in A} dim_z A=\max\limits_{z\in A_{reg}}dim_z A.$$
	
	\n An analytic $A$ is said to be pure $k$-dimentional if $A=\bigcup\limits_{i=1}^m A_i$ is a decomposition of $A$ into irreducible branches $A_i$ and $dim A_i=k$ for all $1\leq i\leq m$. 	
\end{definition}

\n Now, we give the notion of a holomorphic mapping on an analytic set.

\begin{definition}\label{de17} Let $A$ be an analytic subset of a complex manifold $X$. A map $f:\ A \to \mathbb{C}^m$ is said to be a holomorphic map if for each point $a\in A$ there exists a neighborhood $U$ of $a$ in $X$, and a holomorphic mapping $g:\ U \to \mathbb{C}^m$ such that $g=f$ on $A\cap U$.
\end{definition}

Finally, we recall the notion of complex spaces and holomorphic maps between complex spaces. We refer readers to \cite{Ge76,GrRe84,GuRo}. 

\begin{definition}\label{de18} Let $X$ be a topological space. $X$ is said to be a complex space if there exists an atlas of charts consisting of analytic sets such that the transition maps are holomorphic in the sense of holomorphic maps on analytic sets.
\end{definition}

\begin{definition}\label{de19} A point $x$ of a complex space $X$ is said to be regular if there is a neighborhood $U\ni x$ in $X$ such that $U$ is bihomorphically equivalent to an open subset of a complex Euclidean space. The complex dimension of $U$ is called to be the dimension of the set $X$ at its regular $x$, and is denoted by $dim_x X$. The set of regular points of  $X$ is denoted by $X_{reg}$. Every point in the complement $X\setminus X_{reg}=:X_{sing}$ is called to be singular points of $X$. 
\end{definition}

\n The following notion is similar as the notion of irreducible points on an analytic subset in Definition \ref{de15}.	

\begin{definition}\label{de110} A complex space $X$ is called irreducible at $x\in X$ if there is a small neighborhood $U\ni x$ in $X$ such that $U$ is irreducible at $x$. The set of irreducible points of  $X$ is denoted by $X_{irre}$. Every point in the complement $X\setminus X_{irre}=:X_{redu}$ is called to be reducible points of $X$. 
\end{definition}

\n 	

\begin{proposition}\label{pro14} 
	Let $X$ be a complex space. Then 
	
	{\bf i)} $X_{sing}$, $X_{redu}$ are also  complex spaces.
	
	{\bf ii)} $X_{reg}\subset X_{irre}$, $X_{redu}\subset X_{sing}$.
\end{proposition}

\n The following definition is important for the proofs of the main results of the paper.
\begin{definition}\label{def19} Let $X$ be a complex space. $X$ is called to be a strong locally irreducible complex space if it satisfies the following condition: for all $x\in X$ there exists a neighbourhood $U$ of $x$ with the embedding  $j: U\hookrightarrow A\subset\mathbb{C}^n$, $A\subset\mathbb{C}^n$ is a $k$-dimensional analytic subset of $\mathbb{C}^n$. Moreover, let $\triangle'\subset\mathbb{C}^k,\triangle''\subset\mathbb{C}^{n-k}$ be polydisks with the centres at $0'\in\mathbb{C}^k,0''\in \mathbb{C}^{n-k}$ respectively, and with sufficiently small radii $r',r'', r''< r'/C, C>0$ is a constant and let $\Pi:A\bigcap(\triangle^k(0',r')\times\triangle^{n-k}(0'',r''))\rightarrow \triangle^k(0',r')$ be a ramified covering with $p$-sheets then for almost everywhere lines $L\subset \triangle^k(0',r')$, $\Pi^{-1}(L)$ is irreducible in $A$.
\end{definition}

\n Throughout this paper if $X$ is a strong locally irreducible complex space, then we denote it by $\wide{X}$. It is clear that if $X$ is a complex manifold then $X$ is a  strong locally irreducible complex space. In the sequel, if $X=A$ is a $k$-dimensional analytic subset of $\C^n, n\geq k$ then we say that $A$ is a $k$-dimensional strong locally irreducible analytic subset which is understood in the above sense. 

\n Next, we give examples about strong locally irreducible analytic subsets and not strong locally irreducible analytic subsets in $\mathbb{C}^3$.
\begin{proposition}\label{pro15}
	Assume that the analytic subset $A\subset\mathbb{C}^3$ given by
	$$A=\{(z,w,\xi)\in\mathbb{C}^3: \xi^2 -z^k- w^l=0\},$$
	
	\n where $l\geq k\geq 2$. Then $A$ is a strong locally irreducible analytic subset at $\{0\}$ if and only if $k$ is not divisible by 2.
\end{proposition} 

\begin{proof}
	Set $f(z,w,\xi)= \xi^2 - z^k- w^l$. Then $f$ is a holomorphic function on $\mathbb{C}^3$ and $A=\{(z,w,\xi): f(z,w,\xi)= 0\}$, $(0,0,0)\in A$. Set $w=az, a\in\mathbb{C}$. Then $f(z,az,\xi)= \xi^2 - z^k( 1+a^l z^{l-k})$. Assume that $f(z, az,\xi)= (\xi+ g(z))(\xi + h(z))$ where $g(z), h(z)$ are holomorphic functions of $z$. From the decomposition of $f(z,az,\xi)$ it follows that
	\begin{equation*}
	\begin{cases}
	g(z) + h(z) =0\\
	g(z)h(z)= -z^k(1+ a^l z^{l-k})
	\end{cases}
	\end{equation*} 	
	
	\n Hence, it follows that $g(z)=-h(z)$  and we get that $g(z)^2=z^k( 1+ a^l z^{l-k})$. This yields $g(z) = \pm z^{\frac{k}{2}}\sqrt{ 1+ a^l z^{l-k}}$. Therefore, we claim that $A$ is a strong locally irreducible analytic subset at $\{0\}$ if and only if $k$ is not divisible by 2 and $A$ is not a strong locally irreducible analytic subset at $\{0\}$ if and only if $k$ is divisible by 2. The conclusion follows. The same result also holds for analytic subsets $A=\{(z,w,\xi)\in\mathbb{C}^3: \xi^2 + z^k + w^l= 0\}$ in $\mathbb{C}^3$ where $l\geq k\geq 2$.	
	
\end{proof}

\begin{definition}\label{de111} Let $X$, $Y$ be complex spaces. A continue map $f:X\to Y$ is said to be a holomorphic map if for each point $a\in X$ there exists a neighborhood $U$ of $a$ and a neighborhood $V$ of $f(a)$ such that $f:U\to V$ is a holomorphic map in the sense of holomorphic maps on analytic sets.
\end{definition}

\section{Plurisubharmonic functions on complex spaces}\label{section3}
\subsection{ Plurisubharmonic functions on complex spaces}

\n First, we introduce the notion of psh functions on a complex space $X$. Note that the definition of plurisubharmonic functions on complex manifolds is naturally extended from the definition of plurisubharmonic functions on an open subset of $\C^n$. Let $X$ be a complex manifold with $dim_{\mathbb{C}}X=n$ and $\varphi:X\rightarrow [-\infty,+\infty)$ be a given function on $X$. $\varphi$ is called to be a plurisubharmonic function on $X$ if:

$\bullet$ $\varphi$ is upper-semicontinuous on $X$. This means that for all $x_0\in X$,
$$\varlimsup\limits_{X\ni x\to x_0} \varphi(x)\leq \varphi(x_0).$$

$\bullet$ For all $x_0\in X$ there exists a neighbourhood $U$ of $x_0$ and a homeomorphic mapping $\beta:U\cong V\subset\C^n$, where $V$ is an open subset in $\C^n$ such that the function $\varphi\circ \beta^{-1}: V\rightarrow [-\infty,+\infty)$ is a plurisubharmonic function on $V$ which has been introduced and investigated in \cite{BT76,BT82,Kl91,De12}. 

Now we introduce the two definitions of plurisubharmonic functions on a complex space $X$. 

\begin{definition}\label{de21a} 
	
	Let $X$ be a complex space and $\varphi: X\rightarrow [-\infty,+\infty)$. The function $\varphi$ is said to be weakly plurisubharmonic on $X$ if $\varphi$ plurisubharmonic on $X_{reg}$ and 
	\begin{equation}\label{eq3.1}
	\varlimsup\limits_{ x\in X_{reg}, x\to x_0}\varphi(x) = \varphi(x_0),
	\end{equation}
	\n for all $x_0\in X$.
	
	The set of weakly plurisubharmonic functions (resp. negative weakly plurisubharmonic functions) on $X$ is denoted by $PSH_{w} (X)$ (resp. $PSH_{w}^- (X)$).
\end{definition}

\begin{definition}\label{de21b} 
	
	Let $X$ be a complex space and $\varphi: X\rightarrow [-\infty,+\infty)$. The function $\varphi$ is said to be plurisubharmonic on $X$ if:
	
	{\bf i)} $\varphi$ is upper semi-continuous on $X$. This means that for every $x_0\in X$,
	
	$$\varlimsup\limits_{ x\in X, x\to x_0}\varphi(x)\leq \varphi(x_0).$$
	
	{\bf ii)} For all holomorphic mapping $f:\triangle\rightarrow X$, $\varphi\circ f$ is subharmonic on the unit disk $\triangle=\{z\in\mathbb{C}: |z|<1\}$ (possibly identically $-\infty$).
	
	The set of plurisubharmonic functions (resp. negative plurisubharmonic functions) on $X$ is denoted by $PSH (X)$ (resp. $PSH^- (X)$).
\end{definition}
\begin{remark}{\rm
		Note that in the two above definitions, the function $\varphi\equiv-\infty$ on some irreducible branches of a complex space $X$ also is considered to be plurisubharmonic on $X$. Moreover, by Theorem {\bf 5.3.1} in \cite{FN80} the definition of plurisubharmonic functions in \cite{FN80} is equivalent to Definition \ref{de21b}.}
\end{remark}
\begin{example}\label{vd2}{\rm 1) Let $X$ be a complex space and $f\in\mathcal{O}_{X}$ be a holomorphic function on $X$. Then it is easy to see $\varphi=\log|f|$ is a plurisubharmonic function on $X$. Similarly, assume that $f_1,\ldots,f_m$ be holomorphic functions on $X$. Then $\varphi(z)=\log\bigl(|f_1(z)|^2+\cdots+|f_m(z)|^2\bigl)$ also is  a plurisubharmonic function on $X$.
		
		\n 2) Let $A=\{(z_1,z_2)\in \mathbb{C}^2: z_1 z_2=0\}$ be as in Example \ref{vd1}. Then $A$ is a complex space. The set $A_{sing}=\{(0,0)\}$. The function $\varphi(z_1,z_2)=\log|z_1|$ is a plurisubharmonic function on $A$.}
\end{example}

\n Next, we have the following result.
\begin{proposition}\label{pro16}
	Let $X$ be a complex space. Then 
	
	{\bf i)} $PSH (X)\subset PSH_{w}(X)$.
	
	{\bf ii)} $PSH_{w} (X) = PSH (X)$ if $X$ is a complex manifold.
\end{proposition}
\begin{proof}
	
	\n {\bf i)} Let $X$ be a complex space and  $\varphi\in PSH(X)$. We need to prove $\varphi\in PSH_{w}(X)$. Indeed, by the Definition \ref{de21b} the function $\varphi$ is upper-semicontinuous on $X$. Hence, for all $x_0\in X$ we have
	\begin{equation}\label{eq3.2}
	\varlimsup\limits_{ X_{reg}\ni x, x\to x_0}\varphi(x)\leq \varphi(x_0).
	\end{equation}

	\n On the other hand, let $f:\triangle\rightarrow X$ be a holomorphic mapping with $f(0)=x_0$. Then $\varphi\circ f:\triangle \rightarrow [-\infty,+\infty)$ is a subharmonic function on $\triangle$. Hence, we get that
	\begin{equation}\label{eq3.3}
	\varphi(x_0)=\varphi\circ f(0)=\varlimsup_{z\in\triangle, z\to 0}\varphi(f(z))\leq \varlimsup\limits_{x\in X_{reg}, x\to x_0}\varphi(x).
	\end{equation}
	
	\n Coupling \eqref{eq3.3} and \eqref{eq3.2} we get \eqref{eq3.1} and the first desired conclusion of  Definition \ref{de21a} is satisfied. It is enough to check that $\varphi$ is plurisubharmonic on $X_{reg}$. Let $x_0\in X_{reg}$. Then there exists a neighbourhood $U\subset X$ of $x_0$ and a biholomorphic mapping $\beta:U\cong V$ from $U$ onto $V$, where $V$ is an open subset of $\C^n$, $n=dim_{x_0}X$. We must to prove the function $\varphi\circ\beta^{-1}: V\rightarrow [-\infty,+\infty)$ is plurisubharmonic on $V$. Let  $a\in V$, $b\in\C^n$ be arbitrary points. Choose $r>0$ small enough such that $a+rb\lambda\subset V$ for all $\lambda\in\triangle$. We will prove that the function $\triangle\ni\lambda\mapsto \varphi\circ\beta^{-1} (a+rb\lambda)$ is subharmonic on $\triangle$. Consider the holomorphic mapping $h:\triangle\rightarrow V$ such that $h(0)= a$ and $h(\lambda)= a+rb\lambda$. Then the mapping $\beta^{-1}\circ h:\triangle\rightarrow X$ is holomorphic. By the hypothesis the map $\varphi\circ\beta^{-1}\circ h:\triangle\rightarrow [-\infty,+\infty)$ is subharmonic on $\triangle$. Hence, the function  $\triangle\ni\lambda\mapsto \varphi\circ\beta^{-1} (a+rb\lambda)$ is subharmonic on $\triangle$ and the desired conclusion follows. 
	
	\n {\bf ii)}  It is obvious because if $X$ is a complex manifold then $X_{reg}=X$ and hence, the desired conclusion follows from Definition \ref{de21a}.
\end{proof}

\n Now we give an example which indicates that the converse claim of i) in Proposition \ref{pro16} is not true.
\begin{example}\label{vd3}{\rm Consider the complex space $X=\{(z_1,z_2)\in\mathbb{C}^2: z_1.z_2=0\}$ in $\mathbb{C}^2$. Then it is easy to see that $X_{reg}= X\setminus\{(0,0)\}$, $X_{sing}=\{(0,0)\}$. Furthermore, $dim_{a}X=1$ for all $a\in X$. Set $X_1=\{(z_1,0): z_1\in\mathbb{C}\}, X_2=\{(0,z_2): z_2\in \mathbb{C}\}$. Take the function
		$$
		\varphi(z_1,z_2)=
		\begin{cases}
		1   &    \text{if} \  (z_1,0)\in X_1,\\
		0   &    \text{if} \ (0,z_2)\in X_2, z_2\ne 0,\\
		1   &   \text{if} \  (z_1,z_2)=(0,0).
		\end{cases}
		$$
		
		\n It is easy to see that $\varphi$ is a weakly plurisubharmonic function on $X$. Indeed, it is clear that $\varphi$ is plurisubharmonic on $X_{reg}$ and for every $z^0=(z^0_1,z^0_2)\in X$ we have 
		$$\varlimsup\limits_{(z_1,z_2)\in X_{reg}, (z_1,z_2)\to z^0}\varphi(z_1,z_2)= \varphi(z^0)$$
		
		\n However, $\varphi$ is not a psh function on $X$. For the sake of seeking a contradiction we assume that $\varphi$ is a plurisubharmonic function on $X$. Then, by Definition \ref{de21b}, let $f:\triangle\rightarrow X$ be a holomorphic function  given by $f(t)= (0,t), t\in\triangle$ then it follows that $\varphi\circ f$ is subharmonic on $\triangle$. But $\varphi\circ f(t)= 0$ if $t\ne 0$ and $\varphi\circ f (0)=1$ which contradicts to the maximal principle. The proof is complete.}
\end{example}

\n In the sequel, we need the following. 
\begin{definition}\label{dn22}  A one-dimension irreducible complex space $\gamma$ is called a curve germ.
\end{definition}
\begin{proposition}\label{pro200}
	Let $X$ be a complex space and $\varphi : X\to [-\infty ,+\infty )$ be an upper-semicontinuous function. Then the two following statements are equivalent 
	
	{\bf i)} $\varphi$ is plurisubharmonic on $X$.
	
	{\bf ii)} $\varphi |_{\gamma}$ is plurisubharmonic on $\gamma$ for all  curve germ  $\gamma$ in $X$.
\end{proposition}
\begin{proof}
	{\bf i)}$\Longrightarrow${\bf ii)} Let $\gamma\subset X$ be a germ curve. Then $\varphi|_{\gamma}:\gamma\rightarrow [-\infty,+\infty)$ is upper-semicontinous. Let $f:\triangle\rightarrow \gamma$ is a holomorphic mapping. By the hypothesis $\varphi\circ f$ is subharmonic on $\triangle$. Hence, $\varphi|_{\gamma}$ is plurisubharmonic on $\gamma$ and the proof is complete.
	
	\n {\bf ii)}$\Longrightarrow${\bf i)} It remains to check that for all holomorphic mapping $ f:\triangle\rightarrow X$, the function $\varphi\circ f$ is subharmonic on $\triangle$. If $f(t)=constant$ for all $t\in\triangle$ then the conclusion is clear. Let $f(0)=x_0\in X$. Without loss of generality we may assume that $f_{*}(0)$ is invertible (The map $f_{*}$ is defined in \cite{GuRo}, Chapter V, A.11). According to Theorem A.17 in Chapter V in \cite{GuRo} there exists a neighborhood $\triangle(r)=\{z\in\mathbb{C}:|z|<r\}$, $0<r<1$ of $0\in \triangle$ such that $f:\triangle(r)\rightarrow X$ is an immersion. Thus we achieve that $f(\triangle_{r})$ is a curve germ $\gamma$ in $X$. Then by the hypothesis $\varphi\circ f$ is subharmonic on $\triangle_{r}$ and, hence, the desired conclusion follows. 
\end{proof}

\begin{proposition}\label{pro202}
	
	Let $X$ be a complex space. Asume that $\varphi, \psi\in PSH (X)$. Then 
	
	{\bf i)} $a \varphi + b \psi\in PSH (X)$ for all $a,b\geq 0$.
	
	{\bf ii)} $\max (\varphi,\psi)\in PSH (X)$.
\end{proposition}
\begin{proof}
	{\bf i)} It is obvious.
	
	\n {\bf ii)} Since $\varphi\in PSH(X), \psi\in PSH(X)$ then $\varphi,\psi$ are upper-semicontinuous on $X$. Hence, $\max\{\varphi,\psi\}$ also is upper-semicontinuous on $X$. Assume that $f:\triangle\rightarrow X$ is a holomorphic mapping. Then by the hypothesis $\varphi\circ f$ and $\psi\circ f$ are subharmonic on $\triangle$. This yields that $\max\{\varphi\circ f, \psi\circ f\}$ is subharmonic on $\triangle$. However, $\max\{\varphi\circ f, \psi\circ f\}=\max\{\varphi,\psi\}\circ f$ and the desired conclusion follows.
\end{proof}

\begin{proposition}\label{pro203}
	Let $X$ be a complex space and $\{\varphi_j\}_{j\geq 1}$ be a sequence of plurisubharmonic functions on $X$.  Assume that $\varphi_j\searrow\varphi$ on $X$. Then $\varphi$ is plurisubharmonic on $X$. 
\end{proposition}
\begin{proof}
	\n The proof is the same as in the case $X$ is an open subset of $\C^n$ then we omit.
\end{proof}

\begin{proposition}\label{pro204}
	Let $X$ be a complex space and $X=\bigcup\limits_{i\in I} X_i$, where $X_i$ are irreducible branches of $X$. Then 
	$\varphi\in PSH (X)$ if and only if $\varphi\in PSH (X_i)$ for all $i\in I$.
	
\end{proposition}

\begin{proof}

	If $\varphi\in PSH (X)$ then by Definition \ref{de21b} it follows that $\varphi\in PSH (X_i)$ for all $i\in I$. Conversely, we assume that $\varphi\in PSH (X_i)$ for all $i\in I$. We need to prove that $\varphi\in PSH (X)$. First, we prove that $\varphi$ is upper-semicontinuous on $X$. Let $x_0\in X$ be arbitrary. Then it is easy to see that $\varlimsup\limits_{ X\ni x, x\to x_0}\varphi(x)=\sup\limits_{i}\varlimsup\limits_{X_i\ni x, x \to x_0}\varphi(x)\leq \varphi(x_0)$ because $\varphi\in PSH(X_i)$ for all $i\in I$. Next, take a holomorphic mapping $f:\triangle\rightarrow X$ where $\triangle=\{z\in\mathbb{C}: |z|<1\}$ is the unit disc in $\mathbb C$. We have $f^{-1} (X_i)$ are analytic sets in $\triangle$ and $\triangle=\bigcup\limits_{i\in I} f^{-1} (X_i)$. Hence, it follows that there exists $i\in I$ such that $f^{-1} (X_i) = \triangle$. Moreover, since $\varphi\in PSH (X_i)$ then we deduce that $\varphi\circ f$ is subharmonic on $\triangle$.
\end{proof}

\n From now until to the end of this paper we always need an auxiliary result which is the local parametrization theorem II.{\bf 4.19}, p.119 in \cite{De12}. For convenience to readers we recall its contents here and we use it many times in different places of the paper.

\subsection{Local parametrization theorem and auxiliary functions $\varphi_{aver}$, $\varphi_{\max}$ associated with a given function $\varphi$ on a complex space $X$.\\}

\n Let $A$ be a chart of a complex space $X$ at the point $0\in\mathbb {C}^n$ and, moreover, we may assume that $(A,0)$ is irreducible at $0$. Then by local parametrization theorem II.{\bf 4.19} in \cite{De12} there exists the projection  
\begin{equation}\label{eq3.7}
\Pi:A\cap(\triangle^k(0',r')\times\triangle^{n-k}(0'',r''))\rightarrow \triangle^k(0',r')\subset\mathbb{C}^k,
\end{equation}

\noindent which is a ramified covering with $p$-sheets and $S=\Pi ( A_{sing} )$ is an analytic set in $\triangle^k(0',r')$. For each a function $\varphi: A\rightarrow [-\infty,+\infty)$, we set
$$\varphi_{aver} ( \xi )=\frac 1 p \Bigl\{ \sum\varphi (w): \Pi(w)= \xi , w \in A_{reg} \Bigl\},\ \forall \xi\in\triangle^k(0',r')\setminus S.$$
$$\varphi_{max} ( \xi )=\max \{ \varphi (w): \Pi(w)= \xi , w \in A_{reg} \},\ \forall \xi\in\triangle^k(0',r')\setminus S.$$

\n Next, we give the relation between weakly plurisubharmonic functions introduced in this paper and plurisubharmonic functions on a strong locally irreducible complex space $\wide{X}$. Remark that in \cite{De85} Demailly introduced and investigated weak plurisubharmonic functions and plurisubharmonic functions on a complex space $X$ (See Definition 0.1, p.3 in \cite{De85}). He proved that the two these notions are coincident if $X$ is locally irreducible. To prove the above mentioned result, in \cite{De85} Demailly used the method of the resolution of singularities. Here we approach these results by an another method. Before we state and prove some results about plurisubharmonic functions and  weak plurisubharmonic functions on complex spaces, we need the following result.

\begin{lemma}\label{lemma0} Let $\Delta=\Delta(0,1)$ be the unit disk in $\mathbb{C}$ and $f:\Delta\rightarrow\Delta$ be a holomorphic mapping, $f(0)=0$. Assume that $f:\Delta\setminus \{0\}\rightarrow\Delta\setminus\{0\}$ is a covering map with  $p$-sheets. By $\pi_1( \Delta\setminus \{0\} )$ we denote the fundamental group of $\Delta\setminus \{0\}$. Then $f_{\sharp}(\pi_1(\Delta\setminus\{0\}))= p\mathbb{Z}$ where $\mathbb{Z}$ denotes the set of entire numbers.
\end{lemma}	

\begin{proof}
	Indeed, from the hypothesis it follows that $f(t)=t^p g(t)$ where $g(t)$ is a holomorphic function on $\Delta$ with $g(0)\ne 0$. Note that $ \pi_1(\Delta\setminus \{0\})= \mathbb Z$. Next, we have
	$$f_{\sharp}(\pi_1(\Delta\setminus\{0\}))= f_{\sharp}(\mathbb Z)=\mathbb Z f_{\sharp}(1)= p\mathbb Z.$$
	
	\n The desired conclusion follows. 
\end{proof}
\n     
\begin{theorem}\label{theo21c} 
	Let $\wide{X}$ be a strong locally irreducible complex space and  $\varphi: \wide{X}\rightarrow [-\infty,+\infty)$ be a weakly plurisubharmonic function on $\wide{X}$. Then
	
	{\bf i)} $\varlimsup\limits_{ x\in\gamma\backslash\{ x_0 \}, x\to x_0 }\varphi (x) = \varphi (x_0),$ for all curve germ $\gamma$ in $\wide{X}$ at the point $x_0$.
	
	{\bf ii)} $\varphi: \wide{X}\rightarrow [-\infty,+\infty)$ is plurisubharmonic on $\wide{X}$.
\end{theorem}

\begin{proof}
	{\bf i)} First, we consider the case of $x_0\in \wide{X}_{reg}$. Let $\gamma$ be a curve germ at the point $x_0$. Using Example 4.27 in \cite{De12} (also see arguments in \cite{Chir89}, p. 67-68) we can take a local parametrization $f: \triangle(0,1)\rightarrow \gamma$ such that $f(0) = x_0$. Then by the hypothesis $\varphi$ is a weakly plurisubharmonic function on $\wide{X}$ it follows that $\varphi\circ f\in SH (\triangle(0,1))$. Thus this yields that
	$$\varlimsup\limits_{ x\in\gamma\backslash\{ x_0 \}, x\to x_0 }\varphi (x) = \varlimsup\limits_{ z\in\triangle(0,1)\backslash\{ 0 \}, z\to 0 }\varphi\circ f (z) = \varphi\circ f (0) = \varphi (x_0)$$
	Secondly, we consider the case of $x_0\in \wide{X}_{sing}$. By the above passage \eqref{eq3.7}, there exists a ramified covering with $p$-sheets
	$$\Pi:\wide{X}\cap(\triangle^k(0',r')\times\triangle^{n-k}(0'',r''))\rightarrow \triangle^k(0',r')\subset\mathbb{C}^k,$$ 
	Let $S=\Pi ( \wide{X}_{sing} )$. Then $S$ is an analytic set in $\triangle^k(0',r')$, $k=dim \wide{X}$ at $\{0\}$. Let $D$ be the ramification locus of $\Pi$. Then $D$ is an analytic subset of $\triangle^k(0',r')$ and $S\subset D$. Note that
	$$\Pi: ( \wide{X}\cap(\triangle^k(0',r')\times\triangle^{n-k}(0'',r'')) )\backslash\pi^{-1} (D)\rightarrow \triangle^k(0',r')\backslash D,$$ 
	is a covering with $p$-sheets. As in the above passage \eqref{eq3.7}, we construct $\varphi_{aver}$ and $\varphi_{\max}$ associated to $\varphi$ given by
	$$\varphi_{aver} ( \xi )=\frac 1 p \{ \sum\varphi (w): \Pi(w)= \xi ,\ w \in \wide{X}_{reg} \},\ \forall \xi\in\triangle^k(0',r')\setminus D.$$
	$$\varphi_{max} (\xi)=\max \{ \varphi (w): \Pi(w)= \xi ,\ w \in \wide{X}_{reg} \},\ \forall \xi\in\triangle^k(0', r')\setminus D.$$
	
	\n Then Proposition 2.9.26 in \cite{Kl91} implies that $\varphi_{max}$ and $\varphi_{aver}$ are psh functions on $\triangle^k(0',r')\setminus D$. On the other hand, by Theorem 2.9.22 in \cite{Kl91}, $\varphi_{aver}$ and $\varphi_{max}$ can be extended  to two plurisubharmonic functions on $\triangle^k(0',r')$. We will prove that
	$$\varphi_{aver} ( 0 ) = \varphi_{max} ( 0 ) = \varphi (0).$$
	Obviously, we have $\varphi_{aver}\leq \varphi_{max}$  and 
	$$\varphi_{aver} ( 0 ) \leq \varphi_{max} ( 0 ) \leq \varphi (0).$$
	Since 
	$$\varphi (0) = \varlimsup\limits_{ x\in \wide{X}_{reg}, x\to 0}\varphi ( x ) = \varlimsup\limits_{ x\in \wide{X}_{reg}\backslash \pi^{-1} (D), x\to 0}\varphi ( x )\leq \varphi_{max}(0),$$
	we infer that 
	$$\varphi_{max} ( 0 ) = \varphi (0).$$
	\n Next, we choose a line $L=\{ tu:\ t\in\mathbb C \}$, $u\in \triangle^k(0', r')\setminus \{0\}$ of $\triangle^k(0',r')$, such that $L\cap D = \{0\}$. We set $\tau = \Pi^{-1} ( L\cap \triangle^k(0',r'))$ is a curve germ at the point $0$. We have $$\Pi: \tau\backslash\{0\}\rightarrow L\cap \triangle^k(0',r')\backslash\{0\},$$ 
	is a covering with $p$-sheets. We define $g:\Delta\to L\cap \triangle^k(0',r')$ by $g (t) = t^p u$. Then $g:\Delta\backslash\{ 0 \}\to L\cap \triangle^k(0',r')\backslash\{0\}$  is a covering map with $p$-sheets. Moreover, note that $ L\cap \triangle^k(0',r')\backslash\{0\}$ is isomorphic to $\Delta\setminus\{0\}$. Hence, by using Lemma \ref{lemma0} we have $g_{\sharp}(\pi_1(\Delta\setminus\{0\}))= p\mathbb Z$. On the other hand, because $\tau$ is a curve germ at the point $0$ then there exists a local parametrization $h: \Delta\setminus\{0\}\rightarrow \tau\setminus\{0\}$. Then $\Pi\circ h:\Delta\setminus\{0\}\rightarrow L\cap \triangle^k(0',r')\backslash\{0\}$ is a holomorphic mapping. As above, $ L\cap \triangle^k(0',r')\backslash\{0\}\cong \Delta\setminus\{0\}$. Coupling all above arguments we get that
	\begin{align*}
	\Pi_{\sharp}(\pi_1(\tau\setminus\{0\}))= \Pi_{\sharp}(h_{\sharp}(\Delta\setminus\{0\}))
	= (\Pi\circ h)_{\sharp}(\Delta\setminus\{0\})= p\mathbb Z,
	\end{align*}
	
	\n because  $\Pi\circ h: \Delta\setminus\{0\}\rightarrow L\cap \triangle^k(0',r')\backslash\{0\}\cong\Delta\setminus\{0\}$ is a covering map with p-sheets and using Lemma \ref{lemma0}. Now applying  the lifting theorem in \cite{Bre93} (Theorem 4.1, p.143) to $X=\tau\setminus\{0\}$, $Y=L\cap \triangle^k(0',r')\backslash\{0\}$, $p=\Pi$, $f: W\rightarrow Y$ just is the map $g: \Delta\backslash\{ 0 \}\to L\cap \triangle^k(0',r')\backslash\{0\}$ and $g_{\sharp}(\pi_1(\Delta\setminus \{0\}))= p\mathbb Z = \Pi_{\sharp}(\pi_1(\tau\setminus\{0\}))$. Hence, there exists a continuous map $f:\Delta\backslash\{ 0 \}\to\tau\backslash\{ 0 \}$ such that $\Pi\circ f=g$. Since $\Pi\circ f=g$, it follows that $f$ is holomorphic on $\Delta\backslash\{0\}$. Because $f$ is locally bounded then the first removable singularity theorem implies that $f$ is extended to a holomorphic map on $\Delta$. Moreover,
	$$\varphi_{aver}|_{L\cap \triangle^k(z',\delta)} ( \xi ) = \frac 1 p \{ \sum \varphi|_{\gamma} (w): \Pi(w)= \xi ,\ w \in \tau \} = (\varphi|_{ \tau })_{aver}  ( \xi ) ,$$
	$$\varphi_{max}|_{L\cap \triangle^k(z',\delta)} ( \xi ) = \max\{ \varphi|_{ \gamma } (w): \Pi(w)= \xi ,\ w \in \tau \} = (\varphi|_{ \tau })_{max}  ( \xi ),$$
	for all $\xi\in L\cap \triangle^k(0',r')\backslash\{ 0 \}$.
	We have
	$$\varphi_{aver} ( t u ) = \frac {1} {p} \{ \sum \varphi\circ f ( |t|^{\frac 1 p } e^{i \frac { \arg (t) + 2\pi j} {p} } ):\ 0 \leq  j \leq p-1 \},$$
	$$\varphi_{max} ( t u ) = \max \{ \varphi\circ f ( |t|^{\frac 1 p } e^{i \frac { \arg (t) + 2\pi j} {p} } ):\ 0 \leq j \leq p-1 \},$$
	for all $t\in\Delta\backslash\{ 0 \}$.
	On the other hand, since $\varphi\in PSH (\wide{X}_{reg})$ we have $\varphi\circ f\in SH (\Delta\backslash\{ 0 \})$. We can extend $\varphi\circ f$ to a subharmonic function $\psi$ on $\Delta$. Therefore
	$$\begin{aligned}
	\varphi_{aver} ( 0 ) &= \varlimsup\limits_{t\in\Delta\backslash\{0\},\ t\to 0}\varphi_{aver} ( t u )\\ 
	&= \lim\limits_{r\to 0}\frac {1} {2\pi} \int_{0}^{2\pi} \varphi_{aver} ( re^{i\theta} u )\ d\theta
	= \lim\limits_{r\to 0}\frac {1} {2\pi} \int_{0}^{2\pi} \psi ( r^{\frac 1 p }e^{i\theta} )\ d\theta\\
	&=\psi ( 0 )
	= \varlimsup\limits_{t\in\Delta\backslash\{0\},\ t\to 0}\psi (t)\\
	&= \lim\limits_{r\to 0} \max\{ \psi (t):\ |t|=r^{\frac 1 p}\}
	= \lim\limits_{r\to 0} \max\{ \varphi_{max} ( t u ):\ |t|=r\}\\
	&= \varlimsup\limits_{t\in\Delta\backslash\{0\},\ t\to 0}\varphi_{max} ( t u )
	= \varphi_{max} ( 0 )
	\end{aligned}$$
	Therefore
	$\varlimsup\limits_{ x\in\tau\backslash\{ x_0 \}, x\to x_0 }\varphi (x) = \varphi (x_0),$
	for all curve germ $\tau$ in $\wide{X}$ at the point $x_0$ with $\tau\cap \Pi^{-1} (D)={x_0}$.
	This implies that 
	$$\varphi_{aver} ( \xi )=\frac 1 p \{ \sum\varphi (w): \Pi(w)= \xi ,\ w \in \wide{X} \},\ \forall \xi\in\triangle^k(0',r').$$
	$$\varphi_{max} ( \xi )=\max \{ \varphi (w): \Pi(w)= \xi ,\ w \in \wide{X} \},\ \forall \xi\in\triangle^k(0',r').$$
	We set
	$$\varphi_{min} (\xi)=\min \{ \varphi (w): \Pi(w)= \xi ,\ w\in \wide{X} \},\ \forall \xi\in\triangle^k(0',r').$$
	Take a curve germ $\gamma$ in $\wide{X}$ at the point $x_0 = 0$. According to Remmert-Stein's theorem (see \cite{GuRo}), we have $\Pi(\gamma)$ is a curve germ at $0$. Since $\varphi_{aver},\ \varphi_{max}\in PSH (\triangle^k(0',r'))$, we have
	$$\varlimsup\limits_{z\in \Pi(\gamma)\backslash\{ 0 \},\ z\to 0}\varphi_{aver} (z) = \varlimsup\limits_{z\in \Pi(\gamma)\backslash\{ 0 \},\ z\to 0}\varphi_{max} (z) = \varphi (0).$$ 
	This implies that
	$$\varlimsup\limits_{z\in \Pi(\gamma)\backslash\{ 0 \},\ z\to 0}\varphi_{min} (z) \geq \varlimsup\limits_{z\in \Pi(\gamma)\backslash\{ 0 \},\ z\to 0}[ p \varphi_{aver} (z) - (p-1) \varphi_{max} (z) ]\geq \varphi (0).$$
	Hence
	$$\varlimsup\limits_{z\in \Pi(\gamma)\backslash\{ 0 \},\ z\to 0}\varphi_{min} (z) = \varphi (0).$$ 
	This  implies that
	$$\varlimsup\limits_{x\in\gamma\backslash\{ 0 \},\ x\to 0}\varphi (x) \geq \varlimsup\limits_{x\in\gamma\backslash\{ 0 \},\ x\to 0}\varphi_{min} (\Pi (x)) = \varphi (0).$$ 
	\n Hence, 
	$$\varlimsup\limits_{x\in\gamma\backslash\{ 0 \},\ x\to 0}\varphi (x) = \varphi (0).$$ 
	{\bf ii)} By Proposition \ref{pro200} it suffices to show that $\varphi |_{\gamma}$ is subharmonic on $\gamma$ for all curve germ in $\wide{X}$. First, we condider the case where $\gamma$ is a smooth curve germ at the point $0\in \wide{X}$. We choose coordinate such that $\gamma$ is a line. Let $f:\Delta\to\gamma$ be a linear map such that $f(0)=0$. We will prove that 
	$$\varphi\circ f (z)\leq\frac { 1 } {2 \pi r} \int_{\partial\Delta_r} \varphi\circ f (z+w) dV(w),$$
	for all $0 < r < r_0$. 
	Without loss of generality, we may assume that $f(z) = 0 \in \wide{X}_{sing}$. Again, we consider the projection
	$\Pi:\wide{X}\cap(\triangle^k(0',r')\times\triangle^{n-k}(0'',r''))\rightarrow \triangle^k(0',r')\subset\mathbb{C}^k,$
	which is a ramified covering with $p$-sheets and $S=\Pi ( \wide{X}_{sing} )$ is an analytic set in $\triangle^k(0',r')$. Let $D$ be the ramified locus of $\triangle^k(0',r')$ such that 
	$$\Pi: \big[ \wide{X}\cap(\triangle^k(0',r')\times\triangle^{n-k}(0'',r'')) \big]\backslash\pi^{-1} ( D )\rightarrow \triangle^k(0',r')\backslash D,$$  
	is a covering with $p$-sheets and $S\subset D$. We choose a family of lines $L_j\subset\triangle^k(0',r')$ such that $L_j\cap \triangle^k(0',r')\cap D = \{0\}$ and $L_j$ converges to the $L=\Pi ( \gamma )$. We set $\gamma_j = \Pi^{ -1 } (L_j)$. Define $g_j:\Delta\to L_j$ by $g_j ( t ) = t^p u_j, u_j\in \triangle^k(0', r')\setminus\{0\}$. Because $g_j:\Delta\backslash\{ 0 \}\to L_j\cap \triangle^k(0',r')\backslash\{0\}$  is also a covering map with $p$-sheets then by repeating the same arguments as in the proof of {\bf i)} there exist continuous maps $f_j:\Delta\backslash\{ 0 \}\to\gamma_j\backslash\{ 0 \}$ such that $\Pi\circ f_j=g_j$. Since $\Pi\circ f_j=g_j$, as above, we infer that $f_j$ is holomorphic on $\Delta\backslash\{0\}$. This implies that $f_j$ can be extended to a holomorphic map on $\Delta$. 
	Since $\varphi\circ f_j \in SH (\Delta)$, we have
	$$ \varphi\circ f(0)=\varphi\circ f_j (0) = \varphi_{aver} ( f_j (0) )\leq\frac { 1 } {2\pi r} \int_{\partial\Delta_r} \varphi_{aver} (f_j (w)) dV(w),$$
	for all $0 < r < r_0$.
	Moreover, $\Pi\circ f_j$ converges to $\Pi\circ f$ uniformly on every compact set in $\Delta$. Letting $j\to\infty$, we get
	$$\varphi\circ f ( 0 )\leq\frac { 1 } {2\pi r} \int_{\partial\Delta_{r^{\frac 1 p}}} \varphi\circ f ( w ) dV ( w ),$$
	for all $0 < r < r_0$. Hence $\varphi\circ f \in SH (\Delta)$. In general case, by the proof of the first case we have $\varphi |_{\gamma}\in SH ( \gamma_{reg} )$. Moreover, by statement {\bf i)} we infer that $\varphi |_{\gamma}\in SH ( \gamma )$.
\end{proof}

\section{Lelong number of plurisubharmonic functions on complex spaces}\label{section4}

Let $X$ be a complex space and $\varphi\in PSH(X)$. Now, we recall the notion of Lelong number of  $\varphi$ at a point $z\in X$ introduced and investigated in \cite{De85}, p. 45.

\begin{definition}\label{de31} Let $X$ be a complex space and $\varphi\in PSH(X)$. The Lelong number of $\varphi$ at $z\in X$ is defined by
	$$\nu_{\varphi}(z)= \varliminf\limits_{w\to z} \frac {\varphi (w)}{\log ||f(w)-f(z)||},$$
	where $f: U\to A$ is a chart from a neighborhood $U$ of $z$ to an analytic set $A$ in some complex Euclidean space $\mathbb{C}^m$.
\end{definition} 
\begin{proposition}\label{pro32}
	Let $X$ be a complex space and $\varphi$ be a plurisubharmonic function on $X$. Then the definition of Lelong number of $\varphi$ at $z\in X$ does not depend on the charts.
\end{proposition}

\begin{proof}
	Let  $f: U\to A$ and $g: U\to B$ be two charts from a neighborhood $U$ of $z$ where $A$ and $B$ are two analytic sets in two open sets of $\mathbb{C}^k$ and $\mathbb{C}^m$ respectively. By $C_{a} (A)$ we denote the Zariski tangent space of $A$ at $a$. We can choose two manifolds $M$, $N$ in neighborhoods of $f(z)$ and $g(z)$ in $\mathbb{C}^k$ and $\mathbb{C}^m$ be such that $A\subset M$, $B\subset N$ and $T_{ f(z) } M = C_{ f(z) } (A)$, $T_{ g(z) } N = C_{ g(z) } (B)$ where $T_{f(z)} M$ and $T_{g(z)} N$ are the tangent spaces of $M$ at $f(z)$ and $N$ at $g(z)$ respectively. We find two holomorphic maps $F:M\to N$ and $G:N\to M$ such that $F|_A = g\circ f^{-1}$, $G|_B = f\circ g^{-1}$ and $F=G^{-1}$ (see \cite{Wh65}). Hence, it follows that there exist $C_1, C_2>0$ such that
	$C_1 ||f(w)-f(z)||\leq ||g(w)-g(z)||\leq C_2||f(w)-f(z)||$ for all $w\in U$.
	Therefore,
	$$\varliminf\limits_{w\to z} \frac {\varphi (w)}{\log || f(w) - f(z) ||}=\varliminf\limits_{w\to z} \frac {\varphi (w)}{\log || g(w) - g(z) ||},$$
	
	\n and the proof is complete.
\end{proof}

Next, we deal with the following notion.
\begin{definition}\label{de34}
	Let $X, Y$ be  complex spaces and $f:X\rightarrow Y$ be a holomorphic mapping. For each $z\in X, w\in Y$, define 
	$$
	\textup{mult}_{z,w} f=\varliminf\limits_{t\to z}\frac{\log\|h(f(t))-h(w)\|}{\log\|g(t)-g(z)\|},$$

	\n where $h:\Omega_{f(z)}\rightarrow B\subset\mathbb{C}^m$ and $g: V_{z}\rightarrow A\subset \mathbb{C}^n$ are isomorphisms of neighborhoods $\Omega$ of $f(z)$ and $V$ of $z$ onto analytic subsets $B$ and $A$ in $\mathbb{C}^m$ and $\mathbb{C}^n$ respectively.
\end{definition}

\n The following result is an extension of a respective result for holomorphic functions on an open subset in $\C^n$ (see Example {\bf 2.39}, p. 59 in \cite{GuZe}). 
\begin{proposition}\label{theo33}
	Let $X$ be a complex space and $f:X\to \mathbb C$ be a holomorphic function. Then 
	$$\nu_{ \log |f| } (z) = \textup{mult}_{z,0} f,$$
	for all $z\in X$ where $\textup{mult}_{z,0} f$ denotes the vanishing multiplicity of $f$ at a point $z$.
\end{proposition}
\begin{proof}
	The result is directly derived from Definition \ref{de31} and Definition \ref{de34} when we take $h : \mathbb{C}\longrightarrow \mathbb{C}$ is the identity mapping.
\end{proof}

\n The next result is an extension of Theorem {\bf 2.37}, p.58, in \cite{GuZe} from the case of open subsets in $\mathbb{C}^n$ and $\mathbb{C}^m$ to complex spaces. 
\begin{proposition}\label{theo35}
	Let $X$, $Y$ be complex spaces and $f:X\to Y$ be a holomorphic mapping. Assume that $\varphi\in PSH (Y)$. Then 
	$\nu_{ \varphi\circ f } (x) \geq \nu_{ \varphi } ( f(x) ) \textup{mult}_{ x, f(x) } f,$
	for all $x\in X$.
\end{proposition}
\begin{proof}
	It is clear that $\varphi\circ f$ is a plurisubharmonic function on $X$. Let $x\in X$ and $f(x)\in Y$ be arbitrary. Take neighborhoods $V_x$ of $x$ in $X$ and $\Omega_{f(x)}$ of $f(x)$ in $Y$ such that $f(V_x)\subset\Omega_{f(x)}$ and there exist isomorphisms $g: V_{x}\cong A\subset\mathbb{C}^n$ and $h:\Omega_{f(x)}\cong B\subset\mathbb{C}^m$, where $A$ is an analytic subset in $\mathbb{C}^n$ and $B$ is an analytic subset in $\mathbb{C}^m$. Then by Definition \ref{de31} and Definition \ref{de34} we have
	\begin{align*}
	\nu_{\varphi\circ f}(x)&=\varliminf\limits_{t\to x}\frac{(\varphi\circ f)(t)}{\log\|g(t)-g(x)\|}\\
	&=\varliminf\limits_{t \to x}\frac{\varphi(f(t))}{\log \|h(f(t))- h(f(x))\|}\cdot\frac{\log\|h(f(t))-h(f(x))\|}{\log\|g(t)-g(x)\|}\\
	&\geq\varliminf\limits_{t \to x}\frac{\varphi(f(t))}{\log\|h(f(t))-h(f(x))\|}
	\cdot\varliminf\limits_{t \to x}\frac{\log\|h(f(t))-h(f(x))\|}{\log \|g(t)-g(x)\|}\\
	&\geq \nu_{\varphi}(f(x))\cdot \textup{mult}_{x,f(x)}f,
	\end{align*}
	
	\n and the desired conclusion follows.
\end{proof}
\begin{theorem}\label{theo36} Let $X$ be a complex space and $\varphi\in PSH(X)$. Then
	$$\nu_{\varphi}(z) = \inf\bigg\{ \frac { \nu_{\varphi \circ f } (0) } {\textup{ mult}_{0,z} f }\big|\ f:\triangle\to X \ \ \text {is holomorphic with } f(0)=z\bigg\},$$
	for all $z\in X$.
\end{theorem}

\begin{proof}
	By Proposition \ref{theo35} we get that
	$\nu_{\varphi}(z) \leq \inf\big\{ \frac { \nu_{\varphi \circ f } (0)}{ \textup{mult}_{0,z} f}\big|\ f:\triangle\to X\ \ \text {is holomorphic with}\ f(0)=z\big\}.$
	Now, we will show that there exists a holomorphic mapping $f:\triangle\to X$ such that $f(0)=z$, and 
	$$\nu_{\varphi}(z) = \frac { \nu_{\varphi \circ f } (0) } { \textup{mult}_{0,z} f}.$$ 
	First, we prove the above equality in the case where $X$ is a irreducible analytic set at $z$. Assume that $k=dim_z X$. As in the above passage \eqref{eq3.7}, let 
	$$\Pi: X\cap(\triangle^k(z',r')\times\triangle^{n-k}(z'',r''))\rightarrow \triangle^k(z',r')\subset\mathbb{C}^k,$$ 
	be a ramified covering with $p$-sheets, where $z=(z',z'')$ and $S=\Pi ( X_{sing} )$ is an analytic set in $\triangle^k(z',r')$. Set
	$$\psi ( \xi )=\max\{ \varphi (w): \Pi(w)= \xi , w \in X \},\ \forall \xi\in\triangle^k(z',\delta).$$
	Now we choose a line $L=\{ z'+ tu:\ t\in\mathbb C \}$ of $\triangle^k(z',\delta)$, $u\in \triangle^k(z',\delta)\setminus\{0\}$ such that $L\cap S = \{ z' \}$ and 
	$$\nu_{\psi} (z') = \nu_{ \psi|_{ L\cap \triangle^k(z',\delta) } } (z').$$ 
	We have $\Pi :\Pi^{-1} ( L\cap \triangle^k(z',\delta) )\backslash\{ z \} \to L\cap \triangle^k(z',\delta)\backslash\{ z' \}$ is a covering with $p$-sheets. We can write 
	$$\Pi^{-1} ( L\cap \triangle^k(z',\delta) ) = \bigcup\limits_{l=1}^m \gamma_l,$$
	where $\gamma_l$ are locally irreducible branches of $\Pi^{-1} ( L\cap \triangle^k(z',\delta) )$ at the point $z$. We have $\Pi : \gamma_l \backslash\{ z \} \to L\cap \triangle^k(z',\delta)\backslash\{ z' \}$ is a covering with $p_l$-sheets, with $\sum\limits_{1\leq l\leq m} p_l = p$. For each $1\leq l\leq m$, we set 
	$$\psi_l (\xi )=\max\{ \varphi (w): \Pi(w)= \xi , w \in \gamma_l \},\ \forall \xi\in L\cap \triangle^k(z',\delta).$$
	Since $\psi|_{ L\cap \triangle^k(z',\delta) } 
	= \max \{\psi_1, …, \psi_m\},$ then we have
	$$\nu_{ \psi|_{ L\cap \triangle^k(z',\delta) } } (z') = \min \{\nu_{\psi_1} ( z’ ) , …, \nu_{\psi_m} ( z’ ) \}.$$ 
	Therefore, we can find $s\in\{1,…,m\}$ such that
	$$\nu_{ \psi|_{ L\cap \triangle^k(z',\delta) } } (z') = \nu_{\psi_s} ( z’ ).$$
	We consider the map $g:\triangle\to L\cap \triangle^k(z',\delta)$ defined by $g(t) = z' + t^{p_s} u$. Since $g:\triangle\backslash\{ 0 \}\to L\cap \triangle^k(z',\delta)\backslash\{z'\}$ is also a covering with $p_s$-sheets and, by using the same arguments as in the proof of Theorem \ref{theo21c}, we can construct a continuous map $f:\triangle\backslash\{ 0 \}\to \gamma_s \backslash\{ z \}$ such that $\Pi\circ f=g$. Since $\Pi\circ f=g$, it follows that $f$ is holomorphic on $\triangle\backslash\{0\}$ and, hence,  $f$ can be extended to a holomorphic mapping on $\triangle$. We have $f(0)=z$, $\textup{mult}_{0,z} f = \textup{mult}_{0,z'} g=p_s$ and 
	$$\psi_s (z' + \xi u) = \max\{ \varphi\circ f (z'+|\xi|^{ \frac1 p_s }e^{i \frac { \arg (\xi) + 2 p_s j i} {p_s} } u): 0\leq j\leq p_s - 1\}.$$
	This implies that  
	$$\nu_{ \psi_s }( z’ ) = \frac { 1 } { p_s } \nu_{\varphi \circ f } (0) = \frac { \nu_{\varphi \circ f } (0) } { \textup{mult}_{0,z} f}.$$
	Hence, we obtain that
	$$\nu_{ \varphi }(z) = \nu_{ \psi }(z') = \nu_{ \psi|_{ L\cap \triangle^k(z',\delta) } } (z') = \nu_{\psi_s} ( z’ ) = \frac { \nu_{\varphi \circ f } (0) } { \textup{mult}_{0,z} f}.$$
	In the general case, we can write $X=\bigcup\limits_{l=1}^m X_l$ in a neighborhood of $z$, where $X_l$ are locally irreducible branches of $X$ at the point $z$. From the above proof we have
	$$\aligned\nu_{\varphi}(z) &= \min\limits_{1\leq l\leq m} \nu_{ \varphi|_{X_l} } (z)\\ 
	&= \min\limits_{1\leq l\leq m} \inf\{ \frac { \nu_{\varphi \circ f } (0) } { \textup{mult}_{0,z} f }:\ f:\triangle\to X_l \ \ \text { is holomorphic with } f(0)=z\}\\
	&= \inf\{ \frac { \nu_{\varphi \circ f } (0) } { \textup{mult}_{0,z} f }:\ f:\triangle\to X \ \ \text { is holomorphic with } f(0)=z\}.
	\endaligned$$
\end{proof}

\begin{theorem}\label{theo37} Let $X$ be a complex space and $\varphi\in PSH(X)$. Then
	$$\nu_{\varphi}(z) = \inf\{  \nu_{ \varphi |_{\gamma} } (z) :\ \gamma \text{ is a curve germ at the point } z\},$$
	for all $z\in X$.
\end{theorem}

\begin{proof}
	Obviously, we have 
	$\nu_{\varphi}(z) \leq \inf\{  \nu_{ \varphi |_{\gamma} } (z):\ \gamma \text{ is a curve germ at the point } z\}.$
	Now, we prove that there exist a curve germ $\gamma$ at the point $z$ such that 
	$\nu_{\varphi}(z) =  \nu_{ \varphi |_{\gamma} } (z).$
	First, we prove this assertion for the case where $X$ is an irreducible analytic set at $z$. Set $k=dim_z X$. Let 
	$$\Pi: X\cap(\triangle^k(z',\delta)\times\triangle^{n-k}(z'',\eta))\rightarrow \triangle^k(z',\delta)\subset\mathbb{C}^k,$$ 
	be a ramified covering with $p$-sheets, where $z=(z',z'')$ and $S=\Pi ( X_{sing} )$ is an analytic set in $\triangle^k(z',\delta)$. Set
	$$\psi ( \xi )=\max\{ \varphi (w): \Pi(w)= \xi , w \in X \},\ \forall \xi\in\triangle^k(z',\delta).$$
	Now, we choose a line $L=\{ z'+ tu:\ t\in\mathbb C \}$ of $\triangle^k(z',\delta)$ such that $L\cap S = \{ z' \}$ and $\nu_{\psi} (z') = \nu_{ \psi|_{ L\cap \triangle^k(z',\delta) } } (z')$. We set $\gamma = \Pi^{-1} ( L\cap \triangle^k(z',\delta) )$ is a curve germ at the point $z$. We have 
	$$\Pi: \gamma\rightarrow L\cap \triangle^k(z',\delta),$$ 
	is a ramified covering with $p$-sheets. We can write 
	$$\Pi^{-1} ( L\cap \triangle^k(z',\delta) ) = \bigcup\limits_{l=1}^m \gamma_l,$$
	where $\gamma_l$ are locally irreducible branches of $\Pi^{-1} ( L\cap \triangle^k(z',\delta) )$ at the point $z$. For each $1\leq l\leq m$, we set 
	$$\psi_l (\xi )=\max\{ \varphi (w): \Pi(w)= \xi , w \in \gamma_l \},\ \forall \xi\in L\cap \triangle^k(z',\delta).$$
	Since $\psi|_{ L\cap \triangle^k (z',\delta) } 
	= \max \{\psi_1, …, \psi_m\},$ then we have
	$$\nu_{ \psi|_{ L\cap \triangle^k(z',\delta) } } (z') = \min \{\nu_{\psi_1} ( z’ ) , …, \nu_{\psi_m} ( z’ ) \}.$$  
	Therefore, we can find $s\in\{1,…,m\}$ such that
	$$\nu_{ \psi|_{ L\cap \triangle^k(z',\delta) } } (z') = \nu_{\psi_s} ( z’ ).$$
	Thus we infer that                                                                
	$$\nu_{\varphi|_{ \gamma_s }} (z) = \nu_{ \psi_s }
	(z') = \nu_{ \psi|_{ L\cap \triangle^k (z',\delta) } }
	(z') = \nu_{\psi} (z') = \nu_{ \varphi } (z),$$
	
	\noindent and the desired conclusion follows.
	In the general case, we can write $X=\bigcup\limits_{l=1}^m X_l$ in a neighborhood of $z$, where $X_l$ are irreducible branches of $X$ at the point $z$. From the proof just presented  we get that
	$$\aligned\nu_{\varphi}(z) &= \min\limits_{1\leq l\leq m} \nu_{ \varphi|_{X_l} } (z)\\ 
	&= \min\limits_{1\leq l\leq m} \inf\{  \nu_{ \varphi |_{X_l\cap \gamma} } (z) :\ \gamma \text{ is a curve germ in $X_l$ at the point } z\}\\
	&= \inf\{  \nu_{ \varphi |_{\gamma} } (z) :\ \gamma \text{ is a curve germ at the point } z\}
	\endaligned$$
	
	\noindent and the proof is complete.
\end{proof}

\begin{theorem}[\textbf{Theorem A}]\label{theo38}
	Let $X$ be a complex space and $\varphi: X\rightarrow [-\infty,+\infty)$ be a  plurisubharmonic function on $X$. Then for every $x\in X$ the following conclusions hold:
	\vskip0.1cm
	{\bf i)} $\nu_{\varphi_{aver}}(x)\geq \nu_{\varphi}(x) = \nu_{  \varphi_{max}}(x).$
	\vskip0.1cm
	{\bf ii)} If $\wide{X}$ is a strong locally irreducible complex space then for all $x\in\wide{X}$
	\begin{equation}\label{eq4.2}
	\nu_{\varphi_{aver}}(x) = \nu_{\varphi}(x) = \nu_{  \varphi_{max}}(x).
	\end{equation}
	
\end{theorem}
\begin{proof}
	
	\n {\bf i)} Let $x\in X$ be arbitrary. Let \( A \) be a chart of \( X \) at \( x \), and without loss of generality, we assume that \( x = 0 \). First, we prove the following equality:
	\begin{equation}\label{eq4.3}
	\nu_{\varphi}(0) = \nu_{\varphi_{\max}}(0),
	\end{equation}
	
	\noindent	as stated in equation \eqref{eq4.3}.
	
	\n We now assume that \( X = A \) and that \( \varphi \in \text{PSH}(A) \), where \( \varphi \) is a non-positive plurisubharmonic function on \( A \) (i.e., \( \varphi \in \text{PSH}^{-}(A) \)). As in the previous passage \eqref{eq3.7}, consider the map 
	\[
	\Pi: A \cap (\triangle^k(0',r') \times \triangle^{n-k}(0'',r'')) \rightarrow \triangle^k(0',r') \subset \mathbb{C}^k,
	\]
	which is a ramified covering with \( p \)-sheets. For \( x \in A \cap (\triangle^k(0',r') \times \triangle^{n-k}(0'',r'')) \), we write \( x = (x', x'') \), and the set \( S = \Pi(A_{sing}) \) is an analytic subset of \( \triangle^k(0',r') \). Moreover, for each \( x \in A \), we have \( \Pi(x) = x' = (x_1, x_2, \ldots, x_k) \). By the definition we have
	$$\varphi_{\max}(x')=\max\{\varphi(x',x''): (x',x'')\in A\}.$$
	
	\n Furthermore, as demonstrated in the proof of Theorem II.4.19 in \cite{De12}, for any $x=(x',x'')\in A,$ we have $\|x'\|\leq \|x\|\leq C\|x'\|$, $\|x''\|\leq C\|x'\|$. By Definition \ref{de31} we have
	$$\nu_{\varphi}(0)=\lim\limits_{r\to 0}\frac{\max\limits_{\{x\in A: \|x\|=r\}}\varphi(x)}{\log r}.$$
	
	$$\nu_{\varphi_{\max}}(0)=\lim\limits_{r\to 0}\frac{\max\limits_{\{x\in A: \|x'\|=r\}}\varphi_{\max}(x')}{\log r}.$$
	
	\n Note that if $x\in A, \|x\|=r$ and $x'=\Pi(x)$ then
	$\varphi(x)\leq \varphi_{\max}(x').$
	Hence, 
	$$\max\{\varphi(x): x\in A, \|x\|=r\}\leq \max\{\varphi_{\max}(x'):\Pi(x)= x', \|x'\|\leq r\}.$$	
	
	\n Hence, 
	$$\nu_{\varphi}(0)=\lim\limits_{r\to 0}\frac{\max\limits_{\{x\in A: \|x\|=r\}}\varphi(x)}{\log r}\geq \lim\limits_{r\to 0}\frac{\max\limits_{\{x\in A: \|x'\|=r\}}\varphi_{\max}(x')}{\log r}= \nu_{\varphi_{\max}}(0).$$
	Therefore, to derive \eqref{eq4.3}, it suffices to demonstrate that
	\begin{equation}\label{eq4.10}
	\nu_{\varphi_{\max}}(0) \geq \nu_{\varphi}(0).
	\end{equation}	
	\n We define \(\varphi_{\max}(x') = \max\{\varphi(x', x''): (x', x'') \in A\}\). However, based on the previous arguments, we know that \(\|(x', x'')\| \leq (1 + C)\|x'\|\) for \((x', x'') \in A\). This leads to the inequality 
	$$\varphi_{\max}(x') \leq \max\{\varphi(x): x \in A, \|x\| \leq (1 + C)\|x'\|\}.$$ Thus, we have the following inequality:
	
	\begin{equation}\label{eq411}
	\begin{aligned}
	\max\{\varphi_{\max}(x'): \|x'\| = r, x' = \Pi(x), x \in A\} &\leq \max\{\varphi(x): x \in A, \|x\| \leq (1 + C)r\} \\
	&\leq \max\{\varphi(x): x \in A, \|x\| = (1 + C)r\}.
	\end{aligned}
	\end{equation}
	
	\n From \eqref{eq411} we infer that
	\begin{align*}
	\nu_{\varphi_{\max}}(0)=\lim\limits_{r\to 0}\frac{\max\limits_{\{x\in A: \|x'\|=r\}}\varphi_{\max}(x')}{\log r}&\geq\lim\limits_{r\to 0}\frac{\max\limits_{\{x\in A: \|x\|=(1+C)r\}}\varphi(x)}{\log r}\\ 
	&=\nu_{\varphi}(0).
	\end{align*}
	
	\n Thus \eqref{eq4.10} has been proved and the proof of \eqref{eq4.3} is complete.
	
	\n On the other hand, as demonstrated in the proof of Theorem \ref{theo21c}, both $\varphi_{aver}$ and $\varphi_{\max}$ are extended to psh functions on $\triangle^k(0', r')$. We have the inequality
	
	\begin{equation}\label{eq4.12}
	\varphi_{aver}(x') \leq \varphi_{\max}(x'), \quad \text{for all } x' \in \triangle^k(0', r').
	\end{equation}
	
	\noindent Therefore, from \eqref{eq4.12}, and noting that as $r \to 0$, $\log r < 0$, it follows that
	
	\[
	\nu_{\varphi_{aver}}(0) \geq \nu_{\varphi_{\max}}(0) = \nu_{\varphi}(0).
	\]
	\n {\bf ii)} We now demonstrate that if \(\wide{X}\) is a strong locally irreducible complex space, then for all \(x \in \wide{X}\), the following equalities hold:
	\begin{equation}\label{eq4.13}
	\nu_{\varphi_{\text{aver}}}(x) = \nu_{\varphi_{\max}}(x) = \nu_{\varphi}(x).
	\end{equation}
	
	\noindent Without loss of generality, it suffices to prove that equation \(\eqref{eq4.13}\) holds at \(x = 0\). As before, we can assume that \(\varphi \in PSH^{-}(A)\). Furthermore, we again consider that \(\wide{X} = A \subset \mathbb{C}^n\) is a strong locally irreducible analytic subset at \(\{0\}\) in \(\mathbb{C}^n\). Let $\Pi:A\cap(\triangle^k(0',r')\times\triangle^{n-k}(0'',r''))\rightarrow \triangle^k(0',r')\subset\mathbb{C}^k,$
	which is a ramified covering with $p$-sheets, $k=dim_{\{0\}}A$. Assume that $\varphi: A\rightarrow [-\infty,+\infty)$ is a psh function on $A$. As in Theorem \ref{theo21c} we can define $\varphi_{\max}: \triangle^k(0',r')\rightarrow [-\infty,+\infty)$ and $\varphi_{aver}: \triangle^k(0',r')\rightarrow [-\infty,+\infty)$ given by
	$$\varphi_{\max}(x')=\max\{\varphi(x',x''): (x',x'')\in A\}.$$
	$$\varphi_{aver}(x')=\frac{1}{p}\sum\limits_{(x',x'')\in A}\varphi(x',x'').$$
	
	\n Now, let us recall a result by Siu (see \cite{De12}, Theorem III. {\bf 7.13}), which states that for almost every line \( L \subset \mathbb{C}^k \) passing through \( 0 \), we have the following equalities:
	\[
	\nu_{\varphi_{\max}}(0) = \nu_{\varphi_{\max}}|_{L}(0)
	\]
	and
	\[
	\nu_{\varphi_{\text{aver}}}(0) = \nu_{\varphi_{\text{aver}}}|_{L}(0).
	\]	
	\n By the strong  locally irreducibility of $A$ at $\{0\}$ we can choose a line $0\in L\subset \mathbb{C}^k$ such that the three following assertions hold 
	\begin{equation*}
	\begin{cases}
	\text{$\nu_{\varphi_{\max}}(0)=\nu_{\varphi_{\max}}|_{L}(0)$}\\
	\text{$\nu_{\varphi_{aver}}(0)= \nu_{\varphi_{aver}}|_{L}(0)$}\\
	\text {$\Pi^{-1}(L)$ is irreducible.}
	\end{cases}
	\end{equation*}
	
	\n Set $B=\Pi^{-1}(L)\cap A$. Then $B$ is an irreducible curve in $A$ and $\Pi:B\rightarrow L\cap\triangle^k(0',r')$ is an analytic covering with $p$-sheets. Without loss of generality we may consider that $L\cap\triangle^k(0',r')=\triangle(0,1)$, where $\triangle(0,1)$ is the unit disk in $\mathbb{C}$. Put $\Psi= \varphi|_{B}\in PSH^{-}(B)$. Then
	$$\Psi_{\max}=\varphi_{\max}|_{L\cap\triangle^k(0',r')}.$$
	$$\Psi_{aver}=\varphi_{aver}|_{L\cap\triangle^k(0',r')}.$$
	
	\n As in above arguments, we have
	\begin{equation}\label{eq4.18}
	\nu_{\varphi_{\max}}(0)=\nu_{\varphi_{\max}}|_L(0)=\nu_{\Psi_{\max}}(0).
	\end{equation}
	\begin{equation}\label{eq4.19}
	\nu_{\varphi_{aver}}(0)=\nu_{\varphi_{aver}}|_L(0)=\nu_{\Psi_{aver}}(0).
	\end{equation}
	
	\n Since $B$ is an irreducible curve in $A$ then by a result in \cite{De12} about a parametrization of curves (see II.{\bf 4.27}) there exists a bijective holomorphic mapping $f:\triangle \rightarrow B$ such that $\Pi\circ f(t)=t^p, t\in\triangle$. We also note that $\Psi\circ f:\triangle\rightarrow[-\infty,+\infty)$ is a subharmonic function. Next, we prove that
	\begin{equation}\label{eq4.20}
	\nu_{\Psi_{\max}}(0)=\frac{1}{p}\nu_{\Psi\circ f}(0).
	\end{equation}
	
	\n and\begin{equation}\label{eq4.21}
	\nu_{\Psi_{aver}}(0)=\frac{1}{p}\nu_{\Psi\circ f}(0).
	\end{equation}

	\n Assuming that equations \eqref{eq4.20} and \eqref{eq4.21} have been established, it follows that we have the relation:
	
	\begin{equation}\label{eq4.22}
	\nu_{\Psi_{\max}}(0) = \nu_{\Psi_{\text{aver}}}(0).
	\end{equation}
	
	\n Coupling \eqref{eq4.22} with \eqref{eq4.18} and \eqref{eq4.19} we deduce that 
	$\nu_{\varphi_{aver}}(0)= \nu_{\varphi_{\max}}(0).$.
	
	\n Using {\bf i)}, the desired conclusion follows. Hence, it is enough to prove \eqref{eq4.20} and \eqref{eq4.21} hold. First we prove the equality
	$$\nu_{\Psi_{\max}}(0)=\frac{1}{p}\nu_{\Psi\circ f}(0).$$
	
	\n Note that for $x=(x_1,x'')\in B$ there uniquely $t\in\triangle$ such that $f(t)=x$. On the other hand, we have $\Pi(x)=x_1$. Then we have $\Pi\circ f(t)= \Pi(x)=x_1$. But $\Pi\circ f(t)=t^p$. Thus it follows that $t^p=x_1$ Hence,
	$t=\{|x_1|^{1/p} e^{i\frac{arg x_1+2l\Pi}{p}}: l=0,\ldots, p-1\}.$

	\n It follows that
	\begin{align*}
	\Psi_{\max}(x_1) &= \max\{\Psi(x_1, x''): (x_1, x'') \in B\} \\
	&= \max\{\Psi(f(t)): t^p = x_1\} \\
	&= \max\left\{\Psi \circ f\left( |x_1|^{1/p} e^{i \frac{\arg x_1 + 2l\pi}{p}} \right) : l = 0, \ldots, p-1 \right\}
	\end{align*}
	
	\n  Observe that from \(x = f(t)\), it follows that \(x_1 = \Pi(x) = t^p\), and therefore, \(|x_1| = |t|^p\). Hence, we obtain the following expression:
	\begin{equation}\label{eq413}
	\max\limits_{|x_1| = \delta} \Psi_{\max}(x_1) = \max\limits_{|t| = \delta^{1/p}} \Psi \circ f(t).
	\end{equation}

	\n From \eqref{eq413} and by using the definition of the Lelong number we infer that
	$\nu_{\Psi_{\max}}(0)=\frac{1}{p}\nu_{\Psi\circ f}(0),$
	\n and \eqref{eq4.20} is proved. Now we show that \eqref{eq4.21} holds. Indeed, 
	\begin{align*}
	&\frac{1}{2\pi \delta}\int\limits_{|x_1|=\delta}\Psi_{aver}(x_1) d\sigma(x_1)
	=\frac{1}{2\pi}\int_0^{2\pi}\Psi_{aver}(\delta e^{i\theta})d\theta
	=\frac{1}{2\pi p}\sum\limits_{l=0}^{p-1}\int_{0}^{2\pi} \Psi\circ f\Bigl(\delta^{1/p} e^{i\frac{\theta+ 2l\pi}{p}}\Bigl)d\theta\\
	&=\frac{1}{2\pi}\sum\limits_{l=0}^{p-1}\int\limits_{2l \pi/p}^{2(l+1)\pi/p}\Psi\circ f(\delta^{1/p}e^{i\theta})d\theta
	=\frac{1}{2 \pi}\int\limits_0^{2\pi}\Psi\circ f(\delta^{1/p} e^{i\theta})d\theta
	=\frac{1}{2\pi \delta^{1/p}}\int\limits_{|x_1|=\delta^{1/p}}\Psi\circ f(x_1)d\sigma(x_1)
	\end{align*}
	
	\n where $d\sigma(x_1)$ is the Lebesgue measure on the circle $|x_1|=\delta^{1/p}$. On the other hand, if $\Omega$ is an open subset of $\C^n$ and $\varphi\in PSH^{-}(\O)$ and $a\in \O$ then  as in \cite{De12}, we can define the Lelong number of $\varphi$ at $a$ given by
	\begin{equation}\label{eq4.25}
	\nu_{  \varphi}(a)=\nu_{dd^c\varphi}(a)=\lim\limits_{\delta\to 0}\int\limits_{\{\|z-a\|<\delta\}} dd^c\varphi\wedge(dd^c \log\|z-a\|)^{n-1}.
	\end{equation}
	\n From \eqref{eq4.25} it is not difficult to see that the Lelong number of a psh function $\varphi$ at $a\in\Omega$ can be defined by 
	\begin{equation}\label{eq4.26}
	\nu_{\varphi}(a)=\lim\limits_{\delta\to 0}\frac{\int\limits_{\|z-a\|=\delta}\varphi(z)d\sigma(z)}{\log\delta}.
	\end{equation}
	\n where $d\sigma(z)$ denotes the Lebesgue measure of $\{\|z-a\|=\delta\}$. 
	Using \eqref{eq4.26} and the previously obtained result, we derive the following:
	
	\begin{align*}
	\nu_{\Psi_{aver}}(0) &= \lim_{\delta \to 0} \frac{\frac{1}{2\pi \delta} \int_{|x_1| = \delta} \Psi_{aver}(x_1) d\sigma(x_1)}{\log \delta} \\
	&= \lim_{\delta \to 0} \frac{\frac{1}{2\pi \delta^{1/p}} \int_{|x_1| = \delta^{1/p}} \Psi \circ f(x_1) d\sigma(x_1)}{\log \delta} \\
	&= \lim_{\delta \to 0} \frac{\frac{1}{2\pi \delta^{1/p}} \int_{|x_1| = \delta^{1/p}} \Psi \circ f(x_1) d\sigma(x_1)}{p \log \delta^{1/p}} \\
	&= \frac{1}{p} \nu_{\Psi \circ f}(0),
	\end{align*}
	
	\noindent	and hence, we obtain \eqref{eq4.21}. The proof of Theorem \ref{theo38} is now complete.
\end{proof}

\begin{theorem}\label{theo39}
	Let $\wide{X}$ be a strong locally  irreducible complex space. Assume that $\varphi,\psi\in PSH (\wide{X})$. Then 
	
	{\bf i)} $\nu_{ a \varphi } (z) = a \nu_{ \varphi } ( z ),$ for all $a\geq 0$, $z\in \wide{X}$.
	
	{\bf ii)} $\nu_{ \varphi +\psi  }(z) =  \nu_{ \varphi }(z) +  \nu_{\psi } (z),$ for all  $z\in \wide{X}$.
	
	{\bf iii)} $\nu_{ \max ( \varphi , \psi ) } (z) = \min \{ \nu_{ \varphi } ( z ) , \nu_{ \psi } ( z ) \},$ for all $z\in \wide{X}$.
\end{theorem}

\begin{proof} {\bf i)} is trivial.	
	
	\noindent {\bf ii)} We are given a chart \( f: U \to A \), where \( U \) is a neighborhood of a point \( z \in \mathbb{C}^n \) and \( A \) is an analytic set in an open subset of \( \mathbb{C}^n \). The point \( f(z) = a \in A \), and we are further given that \( k = \dim_a A \), the dimension of \( A \) at \( a \).
	Let
	
	\[
	\Pi: A \cap (\triangle^k(a', \delta) \times \triangle^{n-k}(a'', \eta)) \rightarrow \triangle^k(a', \delta) \subset \mathbb{C}^k,
	\]
	
	\noindent be  a ramified covering with \( p \)-sheets where \( a = (a', a'') \) and $S= \Pi(A_{sing})$ is an analytic set in $\triangle(a',\delta)$. 
	Putting
	$\varphi_{aver} ( \xi )= \frac 1 p \{ \sum\varphi\circ f^{-1} (w): \Pi(w)= \xi , w \in A_{reg} \},\ \forall \xi\in\triangle^k(a',\delta)\setminus S.$
	By Proposition {\bf 2.9.26} in \cite{Kl91} it follows that $\varphi_{aver}$ is a psh function on $\triangle^k(a',\delta)$. Since $(\varphi + \psi)_{aver} = \varphi_{aver} + \psi_{aver}$, we get
	$$\nu_{ ( \varphi +\psi )_{aver}  } ( a' ) =  \nu_{ \varphi_{aver} } ( a' ) +  \nu_{ \psi_{aver} } ( a' ).$$
	Now, by using Theorem \ref{theo38} we infer that
	$\nu_{ \varphi +\psi  } (z) =  \nu_{ \varphi } ( z ) +  \nu_{ \psi } ( z ),$
	for all  $z\in \wide{X}$.
	
	\noindent	{\bf iii)} Let $z\in \wide{X}$ and $f: U\to A$  be a chart from a neighborhood $U$ of $z$ onto an analytic set $A$ in an open set $\Omega$ of $\C^n$. We have 
	$$\begin{aligned}
	\nu_{\max(\varphi,\psi)}(z) &=\varliminf\limits_{w\to z}\frac{\max(\varphi ,\psi) (w)}{\log || f(w) - f(z) ||}
	= \varliminf\limits_{w\to z}(-1)\frac{\max(\varphi(w),\psi (w))}{-\log || f(w) - f(z) ||}
	\\
	&= -\varlimsup\limits_{w\to z}\max\left(\frac{\varphi(w)}{-\log || f(w) - f(z) ||}, \frac{\psi(w)}{-\log || f(w) - f(z) ||}\right)\\
	&=-\max\left(\varlimsup\limits_{w\to z}\frac{\varphi(w)}{-\log || f(w) - f(z) ||},\varlimsup\limits_{w\to z}\frac{\psi(w)}{-\log || f(w) - f(z) ||}\right)\\
	&=\min\left(\varliminf\limits_{w\to z}\frac{\varphi(w)}{\log || f(w) - f(z) ||}, \varliminf\limits_{w\to z}\frac{\psi(w)}{\log || f(w) - f(z) ||}\right)
	=\min\{\nu_{  \varphi}(z), \nu_{\psi}(z)\big\},
	\end{aligned} $$
	\n and the desired conclusion follows.
\end{proof}	
\n Next, we give a result about the upper-semicontinuity of Lelong number on complex spaces in the $L^1_{loc}$-topology. Namely, we have the following theorem. 
\begin{theorem}\label{theo310} Let $X$ be a complex space and $\varphi_j,\varphi\in PSH(X)$ and $\varphi_j\rightarrow\varphi$ in $L^1_{loc}(X_{reg})$. Then for all $z\in X$, 
	$\varlimsup\limits_{j\to\infty}\nu_{\varphi_j}(z)\leq \nu_{\varphi}(z).$
\end{theorem}

\begin{proof} 
	Let \( f: U \to A \) be a chart of a neighborhood \( U \) of \( z \), where \( A \) is an analytic set in an open subset of \( \mathbb{C}^n \). Suppose that \( a = f(z) \in A \). First, we prove the theorem in the case $A$ is an irreducible analytic set at $a$. Set $k=dim_a A$. Let 
	$\Pi:A\cap(\triangle^k(a',\delta)\times\triangle^{n-k}(a'',\eta))\rightarrow \triangle^k(a',\delta)\subset\mathbb{C}^k,$ 
	be a ramified covering with $p$-sheets, where $a=(a',a'')$ and $S=\Pi ( A_{sing} )$ is an analytic set in $\triangle^k(a',\delta)$. Set
	$$\psi ( \xi )=\max\{ \varphi\circ f^{-1} (w): \Pi(w)= \xi , w \in A_{reg} \},\ \forall \xi\in\triangle^k(a',\delta)\setminus S.$$
	$$\psi_j ( \xi )=\max\{ \varphi_j\circ f^{-1} (w): \Pi(w)= \xi , w \in A_{reg} \},\ \forall \xi\in\triangle^k(a',\delta)\setminus S.$$
	Since $\Pi$ is a proper holomorphic mapping it follows that if $\varphi_j\rightarrow\varphi$ in $L^1_{loc}(X_{reg})$ then $ \psi_j\rightarrow \psi $ almost everywhere on $\triangle^k(a',\delta)$. Hence, $\psi_j\rightarrow \psi $ in $L^1_{loc}(\triangle^k(a',\delta))$. This implies that $\varlimsup\limits_{j\to\infty}\nu_{\psi_j}(a')\leq \nu_{\psi}(a')$ (refer \cite{De12}). 
	Hence
	$$\varlimsup\limits_{j\to\infty}\nu_{ \varphi_j\circ f^{-1} }(a)\leq \nu_{ \varphi \circ f^{-1} }(a),$$
	\n and we achieve that 
	$\varlimsup\limits_{j\to\infty}\nu_{\varphi_j}(z)\leq \nu_{\varphi}(z),$
	In the general case, we can write $A=\bigcup\limits_{l=1}^m A_l$ in a neighborhood of $a$, where $A_l$ are irreducible branches of $A$. We have
	$$\nu_{ \varphi }(z) = \nu_{ \varphi\circ f^{-1} }(a) = \min\limits_{1\leq l\leq m} \nu_{ \varphi\circ f^{-1}|_{A_l} } (a).$$ 
	$$\nu_{ \varphi_j }(z) = \nu_{ \varphi_j\circ f^{-1} } (a) = \min\limits_{1\leq l\leq m} \nu_{ \varphi_j\circ f^{-1}|_{A_l} } (a).$$ 
	From the proof of the above case, we have that 
	\[
	\varlimsup_{j \to \infty} \nu_{\varphi_j}(z) \leq \nu_{\varphi}(z).
	\]
\end{proof}

\n Let \( A \) be a \( k \)-dimensional analytic set in \( \mathbb{C}^n \). As in \cite{De12}, we can define the closed positive current $[A]$ of bidegree $(n-k,n-k)$ with locally finite masses. Assume that $\varphi$ is a plurisubharmonic function on $A$ in the sense of Definition \ref{de21b}. Then for each $a\in A$ there exists a neighborhood $V$ of $a$ in $\C^n$ such that $\varp$ can be extended to a plurisubharmonic function on $V$ which we still denote it by $\varp$. Thus we can define the closed positive current $T=dd^c(\varp[A])= dd^c \varp\wedge [A]$ of bidegree $(n-k+1,n-k+1)$ on a neighborhood of $a$. By setting
\begin{equation}\label{eq4.27}
\begin{aligned}
\bar{\nu}_{\varp}(a)&=\lim\limits_{r\to 0}\int\limits_{\|z-a\|<r}T\wedge(dd^c\log\|z-a\|)^{k-1}
=\lim\limits_{r\to 0}\int\limits_{\|z-a\|<r} dd^c\varp\wedge(dd^c\log\|z-a\|)^{k-1}\wed [A]\\
&=\int\limits_{\{a\}}dd^c \varp\wedge(dd^c\log\|z-a\|)^{k-1}\wed [A].
\end{aligned}
\end{equation}	
\begin{definition}\label{de35}
	\n $\bar{\nu}_{\varphi}(a)$ is called the projective mass of $\varp$ at $a$.	
\end{definition}
\begin{remark}\label{nx3}{\rm $\bar{\nu}_{\varphi}(a)$ is the Lelong number at $\{a\}$ of the positive closed current $T=dd^c \varp \wedge [A]$.}
\end{remark}
\noindent The following theorem is one of the key results of the paper.
\begin{theorem}[\textbf{Theorem B}]\label{theo311}  Let $A$ be a $k$-dimensional analytic subset in $\C^n, n\geq k$ and $\varphi\in PSH(A)$ be a psh function on $A$. Then the following equality holds for all $a\in A$
	\begin{equation}\label{eq4.28}
	\bar{\nu}_{\varphi}(a)=\textup{mult}(A,a).\nu_{\varphi_{aver}}(a).
	\end{equation}

	\n Furthermore, if $\wide{A}$ is a $k$-dimensional strong locally irreducible analytic subset of $\C^n$ then for all $a\in \wide{A}$,
	$$\bar{\nu}_{\varphi}(a)=\textup{mult}(\wide{A},a).\nu_{\varphi}(a).$$
\end{theorem}	
\begin{proof}
	Without loss of generality we may assume $a=0$. Set $p=\textup{mult}(A,0)$. As in the end of page 213 in \cite{De12}, we have $\textup{mult}(A,0)=\nu([A],0)$. On the other hand, if $A$ is an analytic subset with $dim_{\{0\}}A=k$, $p=\textup{mult}(A,0)$ then by Thie's theorem (see III.{\bf 7.7} in \cite{De12}) there exists a ramified covering with $p$-sheets
	$$\Pi:A\cap(\triangle^k(0',r')\times\triangle^{n-k}(0'',r''))\rightarrow \triangle^k(0',r')\subset\mathbb{C}^k,$$	
	
	\n where $ z= (z',z'')\in A$ satisfying $\|z'\|\leq \|z\|\leq C\|z'\|, \|z''\|\leq C\|z'\|$. Hence, we have $\|z'\|\leq \|z\|\leq (1+C)\|z'\|$. First we prove that \eqref{eq4.28} holds when $\varp\in C^{\infty}(A)$.	By \eqref{eq4.27} we have
	\begin{equation}\label{eq4.30}
	\bar{\nu}_{\varp}(0)=\int\limits_{\{0\}} dd^c\varp\wedge(dd^c\log\|z\|)^{k-1}\wed [A].
	\end{equation}

	\n On the other hand, we have $\lim\limits_{z\to 0}\frac{\log\|z\|}{\log \|z'\|}=1$ then using comparison theorem for Lelong number in \cite{De12} (see III.{\bf 7.1}) we deduce that \eqref{eq4.30} has the form
	\begin{equation}\label{eq4.31}
	\begin{aligned}
	\bar{\nu}_{\varp}(0)&=\int\limits_{\{0\}} dd^c\varp\wedge(dd^c\log\|z\|)^{k-1}\wed [A]
	=\int\limits_{\{0\}} dd^c\varp\wedge (dd^c\log\|z'\|)^{k-1}\wedge [A]\\
	&=\lim\limits_{r \to 0}\frac{1}{r'^{2(k-1)}}\int\limits_{\{\|z'\|<r'\}} dd^c\varp\wedge(dd^c\|z'\|^2)^{k-1}\wedge [A]\\&
	= \lim\limits_{r' \to 0}\frac{1}{r'^{2(k-1)}}\int\limits_{A_{reg}\cap \{\|z'\|<r'\}} dd^c\varp\wedge(dd^c\|z'\|^2)^{k-1},
	\end{aligned}	
	\end{equation}
	
	\n where the third equality is deduced from formula (5.6). p. 202 in \cite{De12} and the fourth equality follows from Definition III.{\bf2.B} in \cite{De12}. However, since the subset \( A_{sing} \cap \{ \|z'\| < r' \} \) has Lebesgue measure equal to zero, we can rewrite equation \eqref{eq4.31} as follows:
	$$\bar{\nu}_{\varp}(0)=\lim\limits_{r \to 0}\frac{1}{r'^{2(k-1)}}\int\limits_{A\cap \{\|z'\|<r'\}} dd^c\varp\wedge(dd^c\|z'\|^2)^{k-1}.$$
	
	\n Next, by the definition of the Lelong number of $\varp_{aver}$ at $\{0\}$ then we have
	$$\nu_{\varp_{aver}}(0)=\lim\limits_{r \to 0}\frac{1}{r'^{2(k-1)}}\int\limits_{ \{\|z'\|<r'\}} dd^c\varp_{aver}\wedge(dd^c\|z'\|^2)^{k-1}.$$
	
	\n Let $S$ be the ramified locus of the ramified covering
	$$\Pi:A\cap(\triangle^k(0',r')\times\triangle^{n-k}(0'',r''))\rightarrow \triangle^k(0',r')\subset\mathbb{C}^k,$$
	
	\n then $S\subset\triangle^k(0',r')$ such that $\Pi:A\cap(\triangle^k(0',r')\times\triangle^{n-k}(0'',r''))\setminus \Pi^{-1}(S)\rightarrow \triangle^k(0',r')\setminus S$ is a covering with $p$-sheets. Write $\Delta^k(0',r')\setminus S =\bigcup\limits_{i=1}^\infty U_i$ , where $U_i$ are open subsets of $\Delta^k(0',r')\setminus S$. Note that $\Pi|_{\Pi^{-1}(U_i)}: \Pi^{-1}(U_i)\rightarrow U_i$ is a covering with $p$-sheets. If set $V_1=U_1, V_2=U_1\setminus U_2,\ldots, V_k=U_k\setminus \bigcup\limits_{l=1}^{k-1} U_l,...$ then we have $\triangle^k(0',r')\setminus S=\bigsqcup\limits_{i=1}^{\infty}V_i, V_s\cap V_t=\emptyset$ when $s\ne t$. Moreover, $\Pi^{-1}(V_i)=\bigcup\limits_{j=1}^p A_{ij}, \Pi|_{A_{ij}} :A_{ij}\simeq V_i$, for all $i\geq 1, 1\leq j\leq p$. Next, we have
	\begin{align*}
	&A\cap(\triangle^k(0',r')\times\triangle^{n-k}(0'',r''))\setminus \Pi^{-1}(S)= \Pi^{-1}(\triangle^k(0',r')\setminus S)
	=\bigcup\limits_{i=1}^{\infty}\Pi^{-1}(V_i).
	\end{align*}
	\n Now we have the following equalities
	
	$$\begin{aligned}
	&\int\limits_{A\cap \{\|z'\|<r'\}\setminus \Pi^{-1}(S)} dd^c\varp\wedge (dd^c\|z'\|^2)^{k-1}\\&=\sum\limits_{i=1}^\infty \int\limits_{\Pi^{-1}(V_i)}dd^c\varp\wedge (dd^c\|z'\|^2)^{k-1}=\sum\limits_{i=1}^\infty\sum\limits_{j=1}^p\int\limits_{A_{ij}}dd^c\varp\wedge (dd^c\|z'\|^2)^{k-1}\\
	&=\sum\limits_{i=1}^\infty\sum\limits_{j=1}^p\int\limits_{V_i} (\Pi|_{A_{ij}})^{-1}_{\ast} \Bigl(dd^c\varp\wedge (dd^c\|z'\|^2)^{k-1}\Bigl)=\sum\limits_{i=1}^\infty\sum\limits_{j=1}^p\int\limits_{V_i} dd^c\varp\circ(\Pi|_{A_{ij}})^{-1} \wedge (dd^c\|z'\|^2)^{k-1}\\
	&=\sum\limits_{i=1}^\infty\int\limits_{V_i}\Bigl(dd^c \sum\limits_{j=1}^p \varp\circ(\Pi|_{A_{ij}})^{-1}\Bigl)\wedge(dd^c \|z'\|^2)^{k-1}= p\sum\limits_{i=1}^\infty\int\limits_{V_i} dd^c \varp_{aver}\wedge(dd^c \|z'\|^2)^{k-1}\\
	&= p\int\limits_{\triangle^k(0',r')\setminus S} dd^c \varp_{aver}\wedge(dd^c \|z'\|^2)^{k-1}= p\int\limits_{\triangle^k(0',r')} dd^c \varp_{aver}\wedge(dd^c \|z'\|^2)^{k-1}.
	\end{aligned}$$

	\n Such as, if $\varphi \in \text{PSH}(A) \cap C^{\infty}(A)$, we obtain that
	\begin{equation}\label{eq4.36}
	\int\limits_{\triangle^k(0',r')} dd^c\varp\wedge (dd^c \|z'\|^2)^{k-1}\wedge [A]=p \int\limits_{ \triangle^k(0',r')} dd^c \varp_{aver}\wedge(dd^c \|z'\|^2)^{k-1}.
	\end{equation}

	\n Now assume that $\varp\in PSH(A),\varp\leq 0$ and $\varp\not\equiv-\infty$. We may assume that $\varphi$ is defined on $\triangle^k(0',r')\times\triangle^{n-k}(0'',r'')$. Take $\var$ to be sufficiently small, with $\var > 0$. Consider $\bar{\triangle}^k(0',r'-\var)\times \bar{\triangle}^{n-k}(0'',r''-\var)$. Choose a decreasing sequence $\{\varp_j\}_{j=1}^{\infty}$ of smooth non-positive plurisubhamonic functions defined on a neighbourhood of $\bar{\triangle}^k(0',r'-\var)\times \bar{\triangle}^{n-k}(0'',r''-\var)$ such that $\varp_j \searrow \varp $ as $j\to\infty$ on  $\bar{\triangle}^k(0',r'-\var)\times \bar{\triangle}^{n-k}(0'',r''-\var)$. Then $(\varp_j)_{aver}\searrow \varp_{aver}$ on $\bar{\triangle}^k(0',r'-\var)$. By \cite{BT82} we have
	\begin{equation}\label{eq4.37}
	dd^c\varp_j\wedge (dd^c\|z'\|^2)^{k-1}\wedge[A]\rightarrow dd^c\varp\wedge(dd^c \|z'\|^2)^{k-1}\wedge [A].
	\end{equation}
	\begin{equation}\label{eq4.38}
	dd^c(\varp_j)_{aver}\wed(dd^c\|z'\|^2)^{k-1} \rightarrow dd^c\varphi_{aver}\wed(dd^c \|z'\|^2)^{k-1}.
	\end{equation}
	\n as $j\to \infty$ in the weak*-topology of currents. On the other hand, for $j\geq 1$, set
	$$H_j=\int\limits_{\triangle^k(0',r'-\var)} dd^c\varp_j\wedge (dd^c \|z'\|^2)^{k-1}\wedge [A].$$
	
	$$K_j=p \int\limits_{\triangle^k(0',r'-\var)} dd^c (\varp_j)_{aver}\wedge(dd^c \|z'\|^2)^{k-1}.$$ 
	
	\n By \eqref{eq4.36} we have $H_j=K_j$ for all $j\geq 1$. However, from the convergence in the weak*-topology of currents in equations \eqref{eq4.37} and \eqref{eq4.38}, we infer that
	\begin{align*}
	& p \int\limits_{\triangle^k(0',r'-\var)} dd^c \varp_{aver}\wedge(dd^c \|z'\|^2)^{k-1}\\
	&\leq \lim\inf\limits_{j\to \infty} K_j
	=\lim\inf\limits_{j\to \infty}H_j\leq \lim\sup\limits_{j\to \infty} H_j
	\leq \int\limits_{\bar{\triangle}^k(0',r'-\var)} dd^c\varp\wed(dd^c \|z'\|^2)^{k-1}\wed [A].	
	\end{align*}	
	\n Letting $\var\to 0$ we get the inequality
	\begin{equation}\label{eq4.40}
	p \int\limits_{ \triangle^k(0',r')} dd^c \varp_{aver}\wedge(dd^c \|z'\|^2)^{k-1}\leq \int\limits_{\triangle^k(0',r')} dd^c\varp\wedge (dd^c \|z'\|^2)^{k-1}\wedge [A].
	\end{equation}
	
	\n By repeating the same arguments as in the above passage, we achieve that
	\begin{equation}\label{eq4.41}
	p \int\limits_{ \triangle^k(0',r')} dd^c \varp_{aver}\wedge(dd^c \|z'\|^2)^{k-1}\geq \int\limits_{\triangle^k(0',r')} dd^c\varp\wedge (dd^c \|z'\|^2)^{k-1}\wedge [A].
	\end{equation}
	
	\n Coupling equations \eqref{eq4.40} and \eqref{eq4.41}, we obtain that
	$$\int\limits_{\triangle^k(0',r')} dd^c\varp\wedge (dd^c \|z'\|^2)^{k-1}\wedge [A]= 
	p \int\limits_{ \triangle^k(0',r')} dd^c \varp_{aver}\wedge(dd^c \|z'\|^2)^{k-1}.$$	
	\n Therefore,
	\begin{equation}\label{eq4.43}
	\int\limits_{A\cap \triangle^k(0',r')} dd^c\varp\wedge (dd^c \|z'\|^2)^{k-1}= p \int\limits_{ \triangle^k(0',r')} dd^c \varp_{aver}\wedge(dd^c \|z'\|^2)^{k-1}.
	\end{equation}

	\n Finally, from \eqref{eq4.43} we obtain that
	\begin{align*}
	\bar{\nu}_{\varp}(0)&=\lim\limits_{r \to 0}\frac{1}{r'^{2(k-1)}}\int\limits_{A\cap \{\|z'\|<r'\}} dd^c\varp\wedge(dd^c\|z'\|^2)^{k-1}\\&=
	p\lim\limits_{r \to 0}\frac{1}{r'^{2(k-1)}}\int\limits_{ \{\|z'\|<r'\}} dd^c\varp_{aver}\wedge(dd^c\|z'\|^2)^{k-1}= p \nu_{\varp_{aver}}(0),
	\end{align*}
	\n and the desired conclusion follows. In the case, $\wide{A}$ is a $k$-dimensional strong locally irreducible analytic subset of $\C^n$ then by Theorem \ref{theo38} we have $\nu_{\varphi_{aver}}(a) = \nu_{\varphi}(a)$ for all $a\in\wide{A}$. Hence, the equality 
	$$\bar{\nu}_{\varphi}(a)=\textup{mult}(\wide{A},a).\nu_{\varphi}(a),$$
	
	\noindent always is done for  $a\in\wide{A}$. The proof of Theorem \ref{theo311} is complete.
\end{proof}

\noindent Below we give counterexamples for {\bf ii)} of Theorem \ref{theo38} and Theorem \ref{theo311}.
\begin{counterexample}\label{phvd1}{\rm
		If we remove the hypothesis of the strong local irreducibility of an analytic subset $A$ in {\bf ii)} of Theorem \ref{theo38} and in Theorem \ref{theo311} then {\bf ii)} of Theorem \ref{theo38} and the equality
		$$\bar{\nu}_{\varphi}(a)= \textup{mult}(A,a).\nu_{\varphi}(a), a\in A,$$
		
		\noindent in Theorem \ref{theo311} do not hold.
		
		\noindent Indeed, consider the analytic subset $A=\{(x,y,z)\in\mathbb{C}^3: x^2 - y^2 + z^3=0\}\subset\mathbb{C}^3$ and $b=(0,0,0)$. By Proposition \ref{pro15}, $A$ is not strong locally irreducible. We have $A_{sing}=(0,0,0)= b$. Consider a psh function $\varphi(x,y, z)=\log\big(|(x+y)^2| + |x-y|+ |z^2|\big)$ in $\mathbb{C}^3$ and its restriction on $A$. Near the origin, the projection $\Pi:\triangle^3(0,\delta)\ni(x,y,z)\longmapsto(x,z)\in\triangle^2(0,\delta)\subset\mathbb{C}^2$ induces a ramified cover with $2$-sheets: $\Pi:A\cap\triangle^3(0,\delta)\rightarrow\triangle^2(0,\delta)$. From the equality $x^2 - y^2 +z^3=0$ it follows that $y=\pm\sqrt{x^2+ z^3}$. On the other hand, by Siu's theorem (\cite{De12}, Theorem III. {\bf 7.13}) for almost everywhere complex line $L\subset\mathbb{C}^3$ we have $\nu_{  \varphi}(0)=\nu\big(\varphi|_{L}\big)(0)$. We have
		$$\varphi_{\max}(x,z)=\max\Bigl[\varphi(x,\pm\sqrt{x^2 +z^3},z)\Bigl].$$
		$$\varphi_{aver}(x,z)=\frac{1}{2}\Bigl[\varphi(x,\sqrt{x^2+z^3},z)+ \varphi(x,-\sqrt{x^2+z^3},z)\Bigl].$$
		
		\n Set $x=a z, a\in \mathbb{C}$, $a\ne 0$. Then it follows that
		$$\varphi_{\max}(az,z)=\max\Bigl[\varphi(az,\pm z\sqrt{a^2+ z)},z)\Bigl].$$
		$$\varphi_{aver}(az,z)=\frac{1}{2}\Bigl[\varphi(az,z\sqrt{a^2+z},z)+ \varphi(az,-z\sqrt{a^2+z},z)\Bigl].$$
		
		\n By the definition of $\varphi$ we have
		$$\varphi(az,z\sqrt{a^2+z},z)= \log\Bigl(\bigl|az+z\sqrt{a^2 +z}\bigl|^2 + \bigl|az-z\sqrt{a^2 +z}\bigl|+ |z|^2\Bigl).$$
		$$\varphi(az,-z\sqrt{a^2+z},z)= \log\Bigl(\bigl|az-z\sqrt{a^2 +z}\bigl|^2 + \bigl|az+z\sqrt{a^2 +z}\bigl|+ |z|^2\Bigl).$$
		
		\n Hence, it is not difficult to see that
		$$\varphi_{max}(az,z)\approx \log|z|+O(1),$$
		
		\n for almost everywhere $a\in\mathbb{C}$, $a\ne 0$. Similarly, we get that
		$$\varphi_{aver}(az,z)\approx \frac{3}{2}\log|z|+ O(1),$$
		
		\n for  almost everywhere $a\in\mathbb{C}$, $a\ne 0$. Hence, $\nu_{ \varphi}(0)= \nu_{\varphi_{\max}}(0)= 1$ and $\nu_{\varphi_{aver}}(0)=\frac{3}{2}.$ Thus we obtain that
		$$\nu_{\varphi_{aver}}(0_{\mathbb{C}^2})=\frac{3}{2}\ne 1=\nu_{\varphi}(0_{\mathbb{C}^3}).$$
		
		\noindent Furthermore, it is not difficult to see that $\textup{mult}(A,0_{\mathbb{C}^3})=2$ (using Corollary at page 122 in \cite{Chir89}). By the equation \ref{eq4.28} we have
		$$\bar{\nu}_{\varphi}(0)= \textup{mult}(A,0_{\mathbb{C}^3}). \nu_{\varphi_{aver}}(0)= 2.\frac{3}{2}=3.$$
		
		\noindent However, 
		$$\textup{mult}(A,0_{\mathbb{C}^3}).\nu_{ \varphi}(0)=2.1=2\ne 3=\bar{\nu}_{\varphi}(0)$$ 
		and, hence, the equality $\bar{\nu}_{\varphi}(a)= \textup{mult}(A,a).\nu_{\varphi}(a)$ does not occur. The fully proven counterexamples.}
\end{counterexample}

\n Finally, we prove a result about a theorem of Siu's type on complex spaces. Namely, we have the following result.

\begin{proposition}\label{pro312} Let $\wide{X}$ be a strong  locally  irreducible complex space and  $\varp\in PSH(\wide{X})$. Then the closure of the set
	$\{x\in \wide{X}: \nu_{\varp}(x) \geq c\}$ is a subvariety  of $\wide{X}$, for all $c\geq 0$.
\end{proposition}

\begin{proof} Let $x\in \wide{X}$. Then there exists a neighbourhood $V_{x}$ of $x$ in $\wide{X}$ such that $V_{x}\simeq A$ is analytically isomorphic to an analytic subset $A$, $dim A=k, A\subset\C^n, k\leq n$. Because the analyticity is a local property then without loss of generality we may assume $\wide{X}= A$ is a $k$-dimensional strong locally irreducible analytic subset of $\C^n$ and $\varphi\in PSH(A)$.
	We need prove that the closure of the set $\mathcal{C}=\{x\in A: \nu_{\varphi}(x)\geq c\}$ is an analytic subset in $A$ for all $c\geq 0$. To get the desired conclusion, it suffices to show that for each $x\in A$ there exists a neighbourhood $U_x$ of $x$ such that $\overline{\mathcal{C}}\cap U_x$ is an analytic subset in $U_x$.  First by Theorem \ref{theo311} we have
	$\bar{\nu}_{\varp}(x)= \textup{mult}(A,x).\nu_{\varphi}(x),$	
	for all $x\in A$. On the other hand,  for each $x\in A$ there exists a neighbourhood $U_x$ of $x$ such that
	$$A\cap U_x=\bigcup\limits_{j=1}^{m}\big\{z\in U_x: \textup{mult}(A,z)=j \big\}.$$
	
	\n Now we have
	\begin{equation}\label{eq4.46}
	\begin{aligned}
	\big\{z\in A\cap U_x: \nu_{\varp}(z)\geq c\big\}&=\bigcup\limits_{j=1}^{m}\Bigl[\big\{z\in A\cap U_x: \nu_{\varphi}(z)\geq c \big\}\bigcap\big\{z\in A\cap U_x:\textup{mult}(A,z)=j\big\}\Bigl]\\
	&= \bigcup\limits_{j=1}^{m}\Bigl[\big\{z\in A\cap U_x : \bar{\nu}_{\varphi}(z)\geq cj \big\}\bigcap\big\{z\in A\cap U_x:\textup{mult}(A,z)=j\big\}\Bigl]
	\end{aligned}	
	\end{equation}
	
	\n 
	Following Remark \ref{nx3}, we have that $\bar{\nu}_{\varphi}(a)$ is the Lelong number at $\{a\}$ of the positive closed current $T=dd^c \varp \wedge [A]$ and by applying Siu’s semi-continuity theorem (see Corollary III.{\bf 8.5} in \cite{De12}), we deduce that $$\Bigl[\Bigl\{z\in A\cap U_x: \bar{\nu}_{\varphi}(z)\geq cj \Bigl\}\bigcap\big\{z\in A\cap U_x:\textup{mult}(A,z)=j\big\}\Bigl]$$ is an analytic subset in $U_x$. On the other hand, 
	\begin{equation}\label{eq4.47}
	\Bigl\{ z\in A\cap U_x:
	\textup{mult}(A,z)=j\Bigl\}= \Bigl\{ z\in A\cap U_x:
	\textup{mult}(A,z)\geq j\Bigl\}\diagdown \Bigl\{ z\in A\cap U_x:
	\textup{mult}(A,z)\geq j+1\Bigl\}.
	\end{equation}
	
	\n By replacing \eqref{eq4.47} into the right-hand side of \eqref{eq4.46} we obtain that
	\begin{align*}
	\mathcal{C}\cap U_x&=\{z\in A\cap U_x: \nu_{\varphi}(z)\geq c\}\\
	&=\bigcup\limits_{j=1}^{m}\Bigl[\Bigl\{z\in A\cap U_x: \bar{\nu}_{\varphi}(z)\geq cj \Bigl\}\bigcap\Bigl\{z\in A\cap U_x:\textup{mult}(A,z)\geq j\Bigl\}\Bigl]\\
	&\diagdown\Bigl[\Bigl\{z\in A\cap U_x: \bar{\nu}_{\varphi}(z)\geq cj \Bigl\}\bigcap\Bigl\{z\in A\cap U_x: \textup{mult}(A,z)\geq j+1\Bigl\}\Bigl].
	\end{align*}
	
	\n For $1\leq j\leq m$, setting
	$$\begin{aligned}
	\mathscr{B}_j&=\Bigl[\Bigl\{z\in A\cap U_x: \bar{\nu}_{\varphi}(z)\geq cj \Bigl\}\bigcap\Bigl\{z\in A\cap U_x:\textup{mult}(A,z)\geq j\Bigl\}\Bigl],\\
	\mathscr{D}_j&=\Bigl[\Bigl\{z\in A\cap U_x: \bar{\nu}_{\varphi}(z)\geq cj \Bigl\}\bigcap\Bigl\{z\in A\cap U_x:\textup{mult}(A,z)\geq j+1\Bigl\}\Bigl],
	\end{aligned}$$
	\n then $\overline{\mathcal{C}}\cap U_x=\bigcup\limits_{j=1}^m \overline{\mathscr{B}_j\setminus\mathscr{D}_j}$. By the above arguments we note that
	$\mathscr{B}_j$ and $\mathscr{D}_j$ are analytic subsets and  $\mathscr{B}_j\diagdown \mathscr{D}_j$ is the difference of the two analytic subsets. Hence, Corollary II.{\bf 5.4} in \cite{De12} implies that the closure of $\overline{\mathscr{B}_j\diagdown \mathscr{D}_j}$ is an analytic subset. As a result, it follows that $\overline{\mathcal{C}}\cap U_x = \bigcup\limits_{j=1}^m \overline{\mathscr{B}_j\diagdown \mathscr{D}_j}$ is an analytic subset in $U_x$.  The proof is complete. 
\end{proof}

%\n{\bf Acknowledgement.} The authors would like to thank the anonymous referees very much for their valuable remarks which led to the improvement of the exposition of this article. 

\end{document}